\documentclass[12pt]{amsart}
\usepackage{graphicx}
\usepackage{amssymb}
\usepackage{amsfonts}

\usepackage{tikz}
\usetikzlibrary{positioning,arrows,matrix,chains,shadows,shapes}
\tikzstyle{Point} = [fill, radius=0.08]

\usepackage{pspicture}
\usepackage{pstricks}
\usepackage[matrix,frame,xdvi]{xypic} 

\newdimen\vcadre\vcadre=0.1cm 
\newdimen\hcadre\hcadre=0.1cm 
\def\GrTeXBox#1{\vbox{\vskip\vcadre\hbox{\hskip\hcadre%
      $#1$%
   \hskip\hcadre}\vskip\vcadre}}
\def\arx#1[#2]{\ifcase#1 \relax \or%
  \ar @{-}[#2]  \or%
  \ar @2{-}[#2] \or%
  \ar @{--}[#2] \or%
  \ar @2{.}[#2] \or%
  \ar @{~}[#2]  \fi}

\def\arbg#1#2{
\newdimen\vcadre\vcadre=0.01cm 
\newdimen\hcadre\hcadre=0.01cm 
\xymatrix@R=0.1cm@C=2mm{
 && {\GrTeXBox{#1}}\arx1[ld]\arx1[rd]\\
& {\GrTeXBox{#2}}\arx1[dl]\arx1[dr] && {\GrTeXBox{z}}\\
 {\GrTeXBox{x}} &&  {\GrTeXBox{y}}\\
}
}
\def\arbd#1#2{
\newdimen\vcadre\vcadre=0.01cm 
\newdimen\hcadre\hcadre=0.01cm 
\xymatrix@R=0.1cm@C=2mm{
 & {\GrTeXBox{#1}}\arx1[ld]\arx1[rd]\\
{\GrTeXBox{x}} & & {\GrTeXBox{#2}}\arx1[ld]\arx1[rd]\\
& {\GrTeXBox{y}}&  & {\GrTeXBox{z}}\\
}
}
\def\arbgnew#1#2{
\newdimen\vcadre\vcadre=0.01cm 
\newdimen\hcadre\hcadre=0.01cm 
\xymatrix@R=0.1cm@C=2mm{
 && {\GrTeXBox{#1}}\arx1[ld]\arx1[rd]\\
& {\GrTeXBox{#2}}\arx1[dl]\arx1[dr] && {\GrTeXBox{B}}\\
 {\GrTeXBox{.}} &&  {\GrTeXBox{A}}\\
}
}

\newtheorem{example}{Example}[section]
\newtheorem{note}[example]{Note}
\newtheorem{theorem}[example]{Theorem}

\newtheorem{corollary}[example]{Corollary}

\newtheorem{definition}[example]{Definition}
\newtheorem{proposition}[example]{Proposition}
\newtheorem{algorithm}[example]{Algorithm}
\newtheorem{lemma}[example]{Lemma}

\def\Proof{\noindent \it Proof -- \rm}
\def\qed{\hspace{3.5mm} \hfill \vbox{\hrule height 3pt depth 2 pt width 2mm}
\bigskip}

\def\sylv{{\rm sylv}}
\def\RC{{\rm RC}}
\def\DC{{\rm DC}}

\def\dT{{\rm DT}}
\def\convW{{*_W}}
\def\convP{{*_P}}
\def\K{{\mathbb K}}
\def\PW{{\rm PW}}
\def\pack{{\rm pack}}
\def\PBT{{\bf PBT}}

\def\FQSym{{\bf FQSym}}
\def\FSym{{\bf FSym}}

\def\PQSym{{\bf PQSym}}
\def\CQSym{{\bf CQSym}}
\def\SQSym{{\bf SQSym}}
\def\WQSym{{\bf WQSym}}
\def\ev{{\rm ev}}
\def\pev{{\rm pev}}

\def\ssh{\Cup}

\def\sconc{\bullet}
\def\Std{{\rm Std}}
\def\std{{\rm std}}
\def\Park{{\rm park}}
\def\mPark{m-{\rm park}}

\def\park{{\bf a}}
\def\a{{\bf a}}

\def\EE{{\sf E}}
\def\<{\langle}
\def\>{\rangle}

\def\Z{\operatorname{\mathbb Z}}

\def\F{{\bf F}}

\def\R{{\bf R}}

\def\G{{\bf G}}

\def\M{{\bf M}}
\def\P{{\bf P}}

\def\SG{{\mathfrak S}}

\def\Sym{{\bf Sym}}

\def\NDPF{{\rm NDPF}}
\def\QSym{{\it QSym}}

\def\Q{{\bf Q}}
\def\NN{{\bf N}}

\def\ch{\operatorname {ch}}

\def\eqhs{\equiv_{\rm hs}}
\def\eqms{\equiv_{\rm ms}}
\def\mPBT{{}^m{\bf PBT}}
\def\mPW{{\rm PW}^{(m)}}
\def\mAn{{}^m{A^n}}
\def\mFQSym{{}^m{\bf FQSym}}
\def\mWQSym{{}^m{\bf WQSym}}
\def\mPQSym{{}^m{\bf PQSym}}
\def\mCQSym{{}^m{\bf CQSym}}
\def\mNCSF{{}^m{\bf NCSF}}
\def\mQSym{{}^m{\bf QSym}}
\def\GD{\rm GD}
\def\mDT{{}^m{\mathcal{DT}}}
\def\mR{{}^m{\mathfrak R}}
\def\PF{{\rm PF}}

\def\MM{{\mathcal M}}

\def\E{{\mathbb E}}

\def\park{{\bf a}}    

\def\shuff#1#2{\mathbin{
\hbox{\vbox{ \hbox{\vrule \hskip#2 \vrule height#1 width 0pt
}%
\hrule}%
\vbox{ \hbox{\vrule \hskip#2 \vrule height#1 width 0pt
\vrule }%
\hrule}%
}}}

\def\shuf{{\mathchoice{\shuff{7pt}{3.5pt}}%
{\shuff{6pt}{3pt}}%
{\shuff{4pt}{2pt}}%
{\shuff{3pt}{1.5pt}}}}%
\def\shuffle{\,\shuf\,}

\def\Tabvrule{\vrule width-0.4pt}       
\def\Tabhrule{\hrule \hrule height-0.4pt} 
\def\Tabstrut{\vrule height2.2ex 
                     depth0.8ex  
                     width0ex    
\relax}

\def\PasCase#1{\omit%
            $\vcenter{\hbox {\vbox to 0.4pt{}}
               \hbox{\makebox[3ex]{\Tabstrut$#1$}}}%
               \Tabvrule$}
\def\PasCasePoint{\PasCase{\cdot}}
\def\DessinCarre#1{%
    \vcenter{\hbox{}\hrule
             \hbox{\vrule\makebox[3ex]{\Tabstrut$#1$}\vrule}\Tabhrule}%
             \Tabvrule}
\def\GenRuban#1{\vcenter{\halign{&$\DessinCarre{##}$\cr#1}}\egroup}

\def\sTabvrule{\vrule width-0.4pt}
\def\sTabhrule{\hrule \hrule height-0.4pt}
\def\sTabstrut{\vrule height1.6ex depth0.6ex width0ex \relax}
\def\sDessinCarre#1{%
    \vcenter{\hbox{}\hrule
             \hbox{\vrule\makebox[2.3ex]%
                  {\sTabstrut$\scriptstyle#1$}\vrule}\sTabhrule}%
             \sTabvrule}
\def\sGenRuban#1{\vcenter{\halign{&$\sDessinCarre{##}$\cr#1}}\egroup}

\def\ruban{%
  \bgroup
  \let\ =\omit
  \let\\=\cr
  \let\x=\times
  \let\.=\PasCasePoint
  \offinterlineskip
  \GenRuban}

\def\sruban{%
  \bgroup
  \let\ =\omit
  \let\x=\times
  \let\\=\cr
  \offinterlineskip
  \sGenRuban}

\newdimen\Squaresize \Squaresize=14pt
\newdimen\Thickness \Thickness=0.5pt

\def\Square#1{\hbox{\vrule width \Thickness
   \vbox to \Squaresize{\hrule height \Thickness\vss
      \hbox to \Squaresize{\hss#1\hss}
   \vss\hrule height\Thickness}
\unskip\vrule width \Thickness}
\kern-\Thickness}

\def\Vsquare#1{\vbox{\Square{$#1$}}\kern-\Thickness}

\def\gf#1#2{\genfrac{}{}{0pt}{}{#1}{#2}}

\def\cerp#1#2{\put(#1,#2){\circle*{0.7}}}
\def\cerg#1#2{\put(#1,#2){\circle*{1}}}

\def\filg#1#2#3{\put(#1,#2){\circle*{.7}\Line(#3,#3)}}
\def\film#1#2#3{\put(#1,#2){\circle*{.7}\Line(0,#3)}}
\def\fild#1#2#3{\put(#1,#2){\circle*{.7}\Line(-#3,#3)}}

\def\arbtun{\begin{picture}(1.5,2)
\cerg{.5}{.5}
\end{picture}}
\def\arbtdeuxa{\begin{picture}(1.1,2.8)
\filg00{1.3} \cerg{1.3}{1.3}
\end{picture}}
\def\arbtdeuxb{\begin{picture}(1.2,3.2)
\film00{1.7} \cerg0{1.7}
\end{picture}}
\def\arbtdeuxc{\begin{picture}(2,2.8)
\fild{1.3}{0}{1.3} \cerg{0}{1.3}
\end{picture}}

\def\arbta{\begin{picture}(2,4.1)
\filg00{1.3} \filg{1.3}{1.3}{1.3} \cerg{2.6}{2.6}
\end{picture}}

\def\arbtb{\begin{picture}(1,3.2)
\filg00{1.7} \film{1.7}0{1.7} \cerg{1.7}{1.7}
\end{picture}}

\def\arbtc{\begin{picture}(2,4.8)
\film00{1.8} \filg0{1.8}{1.5} \cerg{1.5}{3.3}
\end{picture}}

\def\arbtd{\begin{picture}(2,2.8)
\filg00{1.3} \fild{2.6}{0}{1.3} \cerg{1.3}{1.3}
\end{picture}}

\def\arbte{\begin{picture}(2,4.5)
\filg00{1.3} \film{1.3}{1.3}{1.7} \cerg{1.3}{3.0}
\end{picture}}

\def\arbtf{\begin{picture}(3,4.3)
\fild{1.3}0{1.3} \filg{0}{1.3}{1.3} \cerg{1.3}{2.8}
\end{picture}}

\def\arbtg{\begin{picture}(2,3.3)
\fild{1.7}0{1.7} \film00{1.7} \cerg0{1.7}
\end{picture}}

\def\arbth{\begin{picture}(1,4.9)
\film00{1.7} \film0{1.7}{1.7} \cerg{0}{3.4}
\end{picture}}

\def\arbti{\begin{picture}(1,4.3)
\filg{0}0{1.3} \fild{1.3}{1.3}{1.3} \cerg{0}{2.8}
\end{picture}}

\def\arbtj{\begin{picture}(3,4.6)
\fild{1.3}0{1.3} \film{0}{1.3}{1.7} \cerg{0}{3.1}
\end{picture}}

\def\arbtk{\begin{picture}(3,4.6)
\film{1.7}0{1.7} \fild{1.7}{1.7}{1.4} \cerg{0.3}{3.1}
\end{picture}}

\def\arbtl{\begin{picture}(3,4.1)
\fild{2.6}0{1.3} \fild{1.3}{1.3}{1.3} \cerg{0}{2.6}
\end{picture}}

\def\arbqa{\begin{picture}(2,3.3)
\fild{1.7}0{1.7} \film00{1.7} \cerg0{1.7}
\filg{0}{-1.7}{1.7}
\end{picture}}
\def\arbqb{\begin{picture}(3,4.1)
\film{1}{0}{1.7} \film1{1.7}{1.7} \cerg{1}{3.4}
\filg{-0.7}{1.7}{1.7}
\end{picture}}
\def\arbqc{\begin{picture}(3,4.1)
\fild{1.4}{-0.4}{1.7} \filg{-0.3}{1.3}{1.3} \cerg{1}{2.6}
\film{-0.3}{-.4}{1.7}
\end{picture}}
\def\arbqd{\begin{picture}(3,4.1)
\film{-0.7}0{1.7} \filg{-0.7}{1.7}{1.4} \cerg{0.7}{3.1}
\film{-0.7}{-1.7}{1.7}
\end{picture}}
\def\arbqe{\begin{picture}(3,4.1)
\filg{-1.6}0{1.3} \filg{-0.3}{1.3}{1.3} \cerg{1}{2.6}
\film{-1.6}{-1.7}{1.7}
\end{picture}}

\def\arbqf{\begin{picture}(3,4.5)
\film{1}{0}{1.7} \fild{1}{1.7}{1.3} \cerg{-0.3}{3.0}
\filg{-1.6}{1.7}{1.3}
\end{picture}}
\def\arbqg{\begin{picture}(3,4.5)
\filg00{1.3} \film{1.3}{1.3}{1.7} \cerg{1.3}{3.0}
\fild{2.6}0{1.3}
\end{picture}}
\def\arbqh{\begin{picture}(3,4.9)
\film{1}{0}{1.7} \film1{1.7}{1.7} \cerg{1}{3.4}
\filg{-0.3}{-1.3}{1.3}
\end{picture}}
\def\arbqi{\begin{picture}(3,4.6)
\filg{-0.3}0{1.3} \film{1}{1.3}{1.7} \cerg{1}{3.1}
\filg{-1.6}{-1.3}{1.3}
\end{picture}}

\def\arbqj{\begin{picture}(3,3.2)
\filg00{1.7} \film{1.3}0{1.7} \cerg{1.3}{1.7}
\fild{2.6}0{1.7}
\end{picture}}

\begin{document}
\title[Hopf algebras of $m$-Permutations]{Hopf algebras of $m$-permutations,
$(m+1)$-ary trees, and $m$-parking functions}
\author[J.-C.~Novelli and J.-Y.~Thibon]%
{Jean-Christophe Novelli and Jean-Yves Thibon}

\address[]{[Novelli, Thibon] Laboratoire d'informatique Gaspard-Monge\\
Universit\'e Paris-Est Marne-la-Vall\'ee \\
5, Boulevard Descartes \\ Champs-sur-Marne \\
77454 Marne-la-Vall\'ee cedex 2 \\
France}
\email[Jean-Christophe Novelli]{novelli@univ-mlv.fr}
\email[Jean-Yves Thibon]{jyt@univ-mlv.fr} 
\date{\today}


\begin{abstract}
The $m$-Tamari lattice of F. Bergeron is an analogue of the classical Tamari
order defined on objects counted by Fuss-Catalan numbers, such as $m$-Dyck
paths or $(m+1)$-ary trees. On another hand, the Tamari order is related to
the product in the Loday-Ronco Hopf algebra of planar binary trees.
We introduce new combinatorial Hopf algebras based on $(m+1)$-ary trees, whose
structure is described by the $m$-Tamari lattices. 

In the same way as planar binary trees can be interpreted as sylvester classes
of permutations, we obtain $(m+1)$-ary trees as sylvester classes of what we
call $m$-permutations. These objects are no longer in bijection with
decreasing $(m+1)$-ary trees, and a finer congruence, called metasylvester,
allows us to build Hopf algebras based on these decreasing trees. At the
opposite, a coarser congruence, called hyposylvester, leads to Hopf
algebras of graded dimensions $(m+1)^{n-1}$, generalizing noncommutative
symmetric functions and quasi-symmetric functions in a natural way. Finally,
the algebras of packed words and parking functions also admit such
$m$-analogues, and we present their subalgebras and quotients induced by
the various congruences.
\end{abstract}

\maketitle
\footnotesize
\tableofcontents
\normalsize

\section{Introduction}

Among the so-called combinatorial Hopf algebras, the Loday-Ronco algebra of
planar binary trees \cite{LR1} plays a prominent role. Aside from its operadic
interpretation (it is the free dendriform algebra on one generator), it admits
the Connes-Kreimer algebra as a quotient, and shares also many features with
the Hopf algebra of standard Young tableaux, which itself explains the famous
Littlewood-Richardon rule for multiplying Schur functions.

Indeed, both algebras can be defined in similar ways, from two congruences on
words over a totally ordered alphabet, respectively the
sylvester\footnote{Literal translation of the French adjective
{\it sylvestre}, and only an indirect tribute to the great algebraic
combinatorist James Joseph Sylvester.} 
congruence \cite{HNT05} and the plactic congruence \cite{Loth2,PoiReut,NCSF6}.
In both cases, the product  of two basis elements spans an interval of an
order on the relevant combinatorial objects, respectively the Tamari order on
binary trees, and the weak (Melnikov \cite{Mel}) order on standard Young
tableaux.

Recently, the Tamari order has been generalized to an infinite series of
lattices, the $m$-Tamari orders \cite{Ber}, defined on combinatorial objects
counted by Fuss-Catalan numbers, for example $(m+1)$-ary trees. This raises
the question of the existence of a Hopf algebra on $(m+1)$-ary trees, for
which the products of basis elements correspond to $m$-Tamari intervals.

We shall give a positive answer to this question. Actually, there are
``$m$-analogues'' of most classical combinatorial Hopf algebras (permutations,
packed words, parking functions ...), the $(m+1)$-ary trees being obtained
as sylvester classes of $m$-permutations\footnote{The existence of a Hopf algebra
based on $(m+1)$-ary trees has been pointed out in \cite{Ler07}, as a result
of the construction of $m$-dendriform operads. Our constructions appear to
be related to a different operad, which will be investigated
in a separate paper \cite{Nov2}.}, 
as for the Loday-Ronco algebra
obtained for $m=1$. Regarding combinatorial Hopf algebras as generalizations
of the algebra of symmetric functions, we can see that these $m$-analogues
are related to its Adams operations $\psi^m$, raising the variables to the
power $m$. 

These $m$-analogues lead naturally to new congruences, called hyposylvester
and metasylvester, which in turn lead to a wealth of new combinatorial Hopf
algebras.

\bigskip
This article is structured as follows. We first recall the basic constructions
of combinatorial Hopf algebras from congruences and polynomial realizations,
and then establish a criterion for a congruence to be compatible with the Hopf
algebra structure (Section \ref{sec:real}). This criterion, which generalizes
the one of \cite{JNZ}, encompasses all previously known examples.

The new congruences are presented in Section \ref{sec:newcong}.
The metasylvester congruence accounts for the fact that $(m+1)$-ary decreasing
trees are no longer in bijection with $m$-permutations for $m>1$. Thus, one
may expect intermediate algebras of decreasing trees lying between naked
$(m+1)$-ary trees and $m$-permutations. We shall see that metasylvester
classes of $m$-permutations are parametrized by $(m+1)$-ary decreasing trees,
and that sylvester classes are unions of metasylvester classes.

The hyposylvester congruence lies between the sylvester congruence and the
hypoplactic congruence. Its introduction is motivated by the fact that there
are always $2^{n-1}$ hypoplactic classes of $m$-permutations for any $m$,
though one may have expected to see new algebras of graded dimensions
$(m+1)^{n-1}$ in the diagram of homomorphisms. The hyposylvester congruence
does exactly this.  For ordinary permutations ($m=1$), its classes coincide
with the hypoplactic classes, which may explain why it has been overlooked.

The graded dimensions of the various Hopf algebras are computed by means of
explicit formulas for the number of congruence classes of a given evaluation.
These results are of course new for the new congruences, but also for the
sylvester congruence, where it arises as an application of our previous
results on noncommutative Lagrange inversion \cite{NTDup}.

In Section \ref{sec:mperm}, we introduce the Hopf algebras of
$m$-permutations, ${}^m\FQSym$ and its dual, and describe its  quotients and
subalgebras induced by the various congruences.  The algebras of trees,
obtained from the sylvester and metasylvester congruences are discussed in
full detail, the remaining cases being  merely sketched.
It should be noted, however, that the isomorphism of the hypoplactic quotient
(or commutative image) with quasi-symmetric functions is given by the
power-sum plethysm operator (or Adams operator $\psi^m$), and that all this
$m$-combinatorics appears to be related to a noncommutative lift of these
plethysm operations.
We shall also see that the sylvester subalgebra of ${}^m\FQSym$ is isomorphic
to a Hopf subalgebra of the Loday-Ronco algebra, which proves that the product
of two basis elements spans an interval of an ideal of the Tamari order, which
is precisely the $m$-Tamari order.

In Section \ref{sec:mpw}, we introduce Hopf algebras of $m$-packed words,
generalizing $\WQSym$ and its dual. The sylvester quotients and subalgebras
generalizes the free tridendriform algebra on one generator (based on
Schr\"oder trees, counted by the little Schr\"oder numbers). These algebras
are based on $m$-Schr\"oder paths, or on plane trees with no vertex of arity
$m$ or less, counted by generalized Runyon numbers. For the other congruences,
we limit our description to the calculation of the graded dimensions.

Finally, in Section \ref{sec:mpark}, we describe the natural Hopf algebra
structures on $m$-parking functions (in the sense of Bergeron \cite{Ber}),
and investigate subalgebras and quotients.
The hyposylvester quotient is interesting, it is based on $(2m+1)$-ary trees.
Nondecreasing $m$-parking functions, which are in bijection with ordinary
nondecreasing parking functions with repeated letters, are also particularly
interesting. They are involved in the description of the images of the
noncommutative Lagrange series by the adjoints $\phi_m$ of the plethysm
operators, and the dual Hopf algebra ${}^m\CQSym^*$ can be identified as the
plethysm $\psi^m(\CQSym^*)$.
We conclude with the observation that the number of hypoplactic classes of
$m$-parking functions can be obtained by means of $m$-Narayana polynomials,
reciprocals of those giving the number of sylvester classes of $m$-packed
words, and provide some tables of these polynomials and of various numbers of
classes.

\medskip
{\footnotesize
{\it Acknowledgements.-}
This research has been partially supported by
the project CARMA of the French Agence Nationale
de la Recherche.
}

\section{Congruences and polynomial realizations}\label{sec:real}

\subsection{Polynomial realizations}

It is often the case that a combinatorial Hopf algebra admits what we call a
\emph{polynomial realization}. This means that some basis can be faithfully
represented by some (in general, noncommutative) polynomials%
\footnote{More precisely, as elements of an inverse limit of polynomial
algebras in the category of graded rings, as in the case of, \emph{e.g.},
ordinary symmetric functions.},
so that the product rule of the algebra is induced by the ordinary
multiplication of polynomials (or concatenation of words), and the coproduct
is obtained by taking the disjoint union of two copies of the alphabet,
endowed with some extra structure (in general, a total order, but other binary
relations may be involved \cite{FNT}).

The basis elements of a polynomial realization are obtained by summing words
sharing some characteristic. This characteristic may be the result of a
\emph{normalization algorithm}, such as standardization \cite{NCSF6}, packing
\cite{NTtri}, parkization \cite{NTPark}, or an equivalence class for some
congruence on the free monoid. Examples include the plactic congruence
\cite{NCSF6,PoiReut,Loth2}, the hypoplactic congruence \cite{NCSF4,Nov}, the
sylvester congruence \cite{HNT05}, and the stalactic congruence
\cite{NTPark}.
Sometimes, there is an analogue of the Robinson-Schensted correspondence
associated with a congruence, and the $Q$-symbol can play the role of a
normalization algorithm.

\subsection{Basic normalization algorithms}
 
In the sequel, we shall only need a countable totally ordered alphabet
$A$, usually labeled by the positive integers.
We denote by $A^*$ the free monoid generated by $A$.

All algebras will be taken over a field $\K$ of characteristic $0$. The
notation $\K\<A\>$ means the free associative algebra over $A$ when $A$ is
finite, and the inverse limit $\varprojlim_B \K\<B\>$, where $B$ runs over
finite subsets of $A$, when $A$ is infinite.

The \emph{evaluation} $\ev(w)$ of a word $w$ is the sequence whose $i$-th term
is the number of occurrences of the letter $a_i$ in $w$.

The \emph{packed evaluation} $\pev(w)$ of $w$ is the composition obtained by
removing the zeros in $\ev(w)$. 

The \emph{standardized word} $\std(w)$ of a word $w\in A^*$ is the permutation
obtained by iteratively scanning $w$ from left to right, and labelling
$1,2,\ldots$ the occurrences of its smallest letter, then numbering the
occurrences of the next one, and so on. Alternatively, $\sigma=\std(w)^{-1}$
can be characterized as the unique permutation of minimal length such that
$w\sigma$ is a nondecreasing word. For example, $\Std(bbacab)=341625$.

The \emph{packed word} $u=\pack(w)$ of a word $w\in A^*$ where
$A=\{a_1\!<\!a_2\!<\!\cdots\}$ is obtained by the following process.
If $b_1\!<\!b_2\!<\!\cdots\!<\!b_r$ are the letters occurring in $w$, $u$ is the
image of $w$ by the homomorphism $b_i\mapsto a_i$.
A word $u$ is \emph{packed} if $\pack(u)=u$. We denote by $\PW$ the
set of packed words.

Packed words of a given evaluation $I\vDash n$ form a permutation
representation of the symmetric group $\SG_n$ (acting on the right), and also
of the $0$-Hecke algebra $H_n(0)$. The $0$-Hecke algebra, acting by
non-invertible sorting operators, is more convenient, as it allows to
distinguish compositions with the same underlying partition.
The noncommutative characteristic of this action
on words of packed evaluation $I$ (and in particular on packed words,
which are by definition the words whose evaluation is a composition)
 is the product of
noncommutative complete functions $S^I$ \cite{NCSF4}. Its dimension is
therefore the multinomial coefficient $n!S^I(\E)$, where $S_i(\E)$ is by
definition $1/i!$.  Therefore, the characteristic of $\PW$ is
\begin{equation}
\ch (\PW) = \sum_I S^I = (1-S_1-S_2-S_3-\cdots)^{-1}
\end{equation}
which implies immediately that the exponential generating series of $|\PW_n|$
is 
\begin{equation}
\sum_{n\ge 0}|\PW_n|\frac{t^n}{n!}=(1-(e^t-1))^{-1}= (2-e^t)^{-1}. 
\end{equation}
This simple calculation will be generalized several times in the
sequel.

\smallskip

For a word $w$ over a totally ordered alphabet in which each element has a
successor, we defined in \cite{NTPark} a notion of \emph{parkized word}
$\Park(w)$, a parking function which reduces to $\std(w)$ when $w$ is a word
without repetition.

For $w=w_1w_2\cdots w_n$ on $\{1,2,\ldots\}$, we set
\begin{equation}
\label{dw}
d(w):=\min \{i | \#\{w_j\leq i\}<i \}\,.
\end{equation}
If $d(w)=n+1$, then $w$ is a parking function and the algorithm terminates,
returning~$w$. Otherwise, let $w'$ be the word obtained by decrementing all
the values of $w$ greater than $d(w)$. Then $\Park(w):=\Park(w')$. Since
$w'$ is smaller than $w$ in the lexicographic order, the algorithm terminates
and always returns a parking function.

\smallskip
For example, let $w=(3,5,1,1,11,8,8,2)$. Then $d(w)=6$ and the word
$w'=(3,5,1,1,10,7,7,2)$.
Then $d(w')=6$ and $w''=(3,5,1,1,9,6,6,2)$. Finally, $d(w'')=8$ and
$w'''= (3,5,1,1,8,6,6,2)$, which is a parking function.
Thus $\Park(w)=(3,5,1,1,8,6,6,2)$.

For a word $w$ over the alphabet $\{1,2,\ldots\}$, we denote by $w[k]$ the
word obtained by replacing each letter $i$ by the integer $i+k$.
If $u$ and $v$ are two words, with $u$ of length $k$, one defines
the \emph{shifted concatenation}
$u\sconc v = u\cdot (v[k])$
and the \emph{shifted shuffle}
$ u\ssh v= u\shuffle (v[k])$,
where $w_1\shuffle w_2$ is the usual shuffle product defined recursively by
\begin{itemize}
\item $w_1 \shuffle \epsilon  = w_1$,\quad $\epsilon \shuffle w_2  = w_2$,
\item $au \shuffle bv = a(u \shuffle bv) + b(au \shuffle v)$,
\end{itemize}
where $w_1=a\cdot u$ and $w_2=b\cdot v$, and both $a$ and $b$ are letters and
$\cdot$ means concatenation.

For example,
\begin{equation}
12\ssh 21 = 12\shuffle 43 = 1243 + 1423 + 1432 + 4123 + 4132 + 4312\,.
\end{equation}

\subsection{Combinatorial Hopf algebras from elementary characteristics}

\subsubsection{Permutations}

Taking the standardization as characteristic, we obtain the algebra
$\FQSym(A)$ of Free Quasi-Symmetric functions, an algebra
spanned by the noncommutative polynomials \cite{NCSF6}
\begin{equation}
\G_\sigma(A)  := \sum_{\gf{w\in A^n}{\std(w)=\sigma}} w,
\end{equation}
where $\sigma$ is a permutation in the symmetric group $\SG_n$.
The multiplication rule is
\begin{equation}
\G_\alpha \G_\beta = \sum_{\gamma\in \alpha * \beta} \G_\gamma,
\end{equation}
where the convolution $\alpha*\beta$ of $\alpha\in\SG_k$ and $\beta\in\SG_l$
is the sum in the group algebra of $\SG_{k+l}$~\cite{MR}
\begin{equation}
\alpha * \beta =
\sum_{\gf{\gamma=uv}{\std(u)=\alpha ;\, \std(v)=\beta}} \gamma\,.
\end{equation}
$\FQSym$ is therefore a polynomial realization of the Malvenuto-Reutenauer
Hopf algebra~\cite{MR}.

Let $\F_\sigma=\G_{\sigma^{-1}}$. 
If one denotes by $A'$ and $A''$ two mutually commuting alphabets
isomorphic to $A$ as ordered sets, and by $A' + A''$ the ordinal sum, 
the coproduct is defined by the ordinal sum of alphabets
\begin{equation}
 \Delta\F_\sigma =
 \F_\sigma(A'+A'')=
 \sum_{u\cdot v=\sigma}\F_{\std(u)}\otimes \F_{\std(v)}\,,
\end{equation}
(where $\cdot$ means concatenation)
under the identification $U(A)\otimes V(A)=U(A')V(A'')$.

A scalar product is defined by
\begin{equation}
\<\F_\sigma\,,\,\G_\tau\>=\delta_{\sigma,\tau}\,,
\end{equation}
and one then has, for all $F,G,H\in\FQSym$ 
\begin{equation}
\<FG,H \>=\<F\otimes G,\Delta H\>\,
\end{equation}
that is, $\FQSym$ is self-dual.
The product formula in the $\F_\sigma$ basis is
\begin{equation}
  \F_\alpha \F_\beta = \sum_{\sigma \in \alpha\ssh\beta} \F_\sigma.
\end{equation}

\subsubsection{Packed words}

Taking the packing algorithm as characteristic, we define polynomials
\begin{equation}
\M_u :=\sum_{\pack(w)=u}w\,,
\end{equation}
for all packed word $u$.
Under the abelianization map $\chi:\ \K\langle A\rangle\rightarrow\K[X]$,
sending the noncommuting variables $a_i$ to mutually commuting variables
$x_i$, the $\M_u$ are mapped to the monomial quasi-symmetric functions $M_I$
($I=(|u|_a)_{a\in A}$ being the evaluation vector of $u$).

The $\M_u$ span a subalgebra of $\K\langle A\rangle$, called $\WQSym$ for Word
Quasi-Symmetric functions~\cite{Hiv}, consisting in the invariants of the
noncommutative version \cite{NCSF7} of Hivert's quasi-symmetrizing action
\cite{Hiv-adv}, which is defined by $\sigma\cdot w = w'$ where $w'$ is such
that $\std(w')=\std(w)$ and $\chi(w')=\sigma\cdot\chi(w)$.
Hence, two words are in the same $\SG(A)$-orbit iff they have the same packed
word.

The product of $\WQSym$ is given by
\begin{equation} 
\label{prodG-wq}
\M_{u'} \M_{u''} = \sum_{u \in u'\convW u''} \M_u\,,
\end{equation}
where the \emph{convolution} $u'\convW u''$ of two packed words
is defined as
\begin{equation} 
u'\convW u''
  = \sum_{\gf{u=v\cdot w}{\pack(v)=u',\,\pack(w)=u''}}
    u\,.
\end{equation}

Evaluating $\M_u$ on an ordinal sum, we get
\begin{equation}
\M_u(A' + A'')
  = \sum_{0\leq k\leq \max(u)}
          \M_{(u|_{[1,k]})}(A') \M_{\pack(u|_{[k+1,\max(u)]})}(A''),
\end{equation}
where $u|_{B}$ denotes the subword obtained by restricting $u$ to the
subset $B$ of the alphabet,
so that the coproduct $\Delta$ is given by
\begin{equation}
\label{coprG-wq}
\Delta \M_u(A)
=\sum_{0\leq k\leq \max(u)}
\M_{(u|_{[1,k]})}\otimes \M_{\pack(u|_{[k+1,\max(u)]})}.
\end{equation}

\subsubsection{Parking functions}
\label{sec-park}

Let $A=\{a_1<a_2<\dots\}$. Taking the parkization algorithm as characteristic,
we can define polynomials
\begin{equation}
\label{realG}
\G_\park(A) = \sum_{\Park(w)=\park} w,
\end{equation}
$\park$ being a parking function.

The $\G_{\park}(A)$ span a subalgebra of $\K\langle A\rangle$.
This algebra is denoted by $\PQSym^*$, and its product rule is given by
\begin{equation}
\label{prodG}
\G_{\park'} \G_{\park''}  
= \sum_{\park \in \park'\convP\park''} \G_\park\,,
\end{equation}
where the \emph{convolution} $\park'\convP\park''$ of two parking functions
is defined as
\begin{equation}
\label{convolP}
 \park'\convP\park'' =
 \sum_{\gf{\park=u\cdot v}{\Park(u)=\park',\,\Park(v)=\park''}}
      \park\,.
\end{equation}

The coproduct is again given by the ordinal sum
\begin{equation}
\label{coprodreal}
\Delta(\G)(A)=\G(A' + A'').
\end{equation}

The dual basis of $\G_\park$ is denoted by $\F_\park$. The product
of $\PQSym$ is a shifted shuffle
\begin{equation}
\label{prodF}
\F_{\park'}\F_{\park''}
  := \sum_{\park \in \park'\ssh\park''}\F_\park\,.
\end{equation}

\subsection{Hopf algebras from congruences}
\label{sub-cong}

Further examples of combinatorial Hopf algebras can be constructed
by quotienting polynomial realizations by various congruences
on $A^*$, or dually, by summing permutations, packed words or parking
functions over congruence classes. All known examples are special cases
of the following construction.

\def\slha{{shuffle-like Hopf algebra }}
\def\slhas{{shuffle-like Hopf algebras }}

Let us define a \emph{\slha} $H$ as a graded connected Hopf algebra, with a
basis $b_w$ indexed by some class of words over the alphabet $A$ of positive
integers such that the product in this basis is given by the shifted shuffle
and the coproduct is of the form
\begin{equation}
\Delta b_w = \sum_{w=w_1\cdot w_2} b_{a(w_1)} \otimes b_{a(w_2)},
\end{equation}
where $a$ is an algorithm such that $a(u\cdot\sigma)=a(u)\cdot\sigma$ for any
permutation $\sigma$ and such that 
$\pack(a(v))=\pack(v)$. 

For example, $\FQSym$, $\WQSym$, and $\PQSym$ 
are \slhas\, since their product in the $\F$ bases is given by the shifted
shuffle and the algorithm used in their coproduct is $\pack$ for the
first two algebras and $\Park$ for the last one.
We shall see that all $m$-generalizations of these algebras also are \slhas.

The following result generalizes and refines~\cite{JNZ}[Chap. 2.1].

\begin{theorem}
\label{th:cong}
Let $H$ be a \slha.
Let $\equiv$ be a congruence on $A^*$ such that
\begin{itemize}
\item two words are equivalent iff they have same evaluation and their packed
words are equivalent,
\item they are such that
\begin{equation}
u\equiv v \implies u|_I\equiv v|_I,
\end{equation}
for any interval $I$ of $A$, the notation $w|_I$ meaning the restriction of
$w$ to $I$, obtained by erasing the letters not in $I$.
\end{itemize}

Then, the elements
\begin{equation}
P_u = \sum_{w\equiv u}b_w
\end{equation}
span a Hopf subalgebra $A$ of $H$.

Moreover, its dual Hopf algebra $A^*$ can be identified with the quotient of
$H^*$ by the relations $b^*_w=b^*_{w'}$ iff $w\equiv w'$.
\end{theorem}

\Proof
Let us first consider a product $P_{u_1} P_{u_2}$ and an element $b_w$
obtained by expanding this product on the $b$ basis. Note that the $P_u$ are
sums over disjoint sets of $b_w$, so that any $b_w$ occurs in exactly one
$P_u$.
Let us now consider $w'\equiv w$. Let $k$ be the size of $u_1$, $\ell$ the
size of $u_2$. Then the restrictions of both $w$ and $w'$ to
$I_1=[1,k]$ are congruent, and the same holds for the restrictions to
$I_2=[k+1,k+\ell]$. So $w|_{I_1}$ and $w'|_{I_1}$ are both congruent to $u_1$
and $w|_{I_2}$ and $w'|_{I_2}$ are both congruent to $u_2[k]$ since
$w|_{I_2}=w_2[k]$. So $b_{w'}$ is also obtained in the expansion of $P_{u_1}
P_{u_2}$, so that the product is indeed a sum of terms $P_u$.

\medskip
Let us now consider the coproduct of a $P_u$ and an element
$b_{w'}\otimes b_{w''}$ appearing in it. Consider the set $S$ of words $w$
with $w\equiv u$ and such that the same
$b_{w'}\otimes b_{w''}$ appears in the coproduct of $b_w$.
Write for each such word $w=w_1\cdot w_2$, so that $a(w_1)=w'$ and
$a(w_2)=w''$.

Consider now another element
$b_{v'}\otimes b_{v''}$ with $v'\equiv w'$ and $v''\equiv w''$. We want to
prove that the coefficient of this element in the coproduct of $P_u$ is the
same as the coefficient of $b_{w'}\otimes b_{w''}$ is the same coproduct.

Since $\ev(v')=\ev(w')$ and $\ev(v'')=\ev(w'')$, one defines two
permutations as the permutations with smallest number of inversions so that
$v'=w'\cdot\sigma_1$ and $v''=w''\cdot \sigma_2$.

Then let $S'=S\cdot\sigma_1\sconc\sigma_2$.
The prefixes $v_1$ of the elements of $S'$ 
satisfy $v_1=w_1\cdot\sigma_1$ and also satisfy $a(v_1)=v'$
by the properties of the algorithm $a$.
So $\pack(v_1)=\pack(v')\equiv\pack(w')=\pack(w_1)$, and since
$\ev(v_1)=\ev(w_1)$, we conclude that $v_1\equiv w_1$. The same holds for the
suffixes, so that all the $v$ also belong to the same equivalence class.

\bigskip
As for the dual Hopf algebra, if one has an injective Hopf algebra map
from $A$ to $S$, then its transpose is a surjective Hopf algebra map from
$S^*$ to $A^*$, hence realizing the dual $A^*$ as claimed. Indeed,  the
$P_u$ are disjoint sums of $b_w$,  so that the first matrix has one $1$
in each column, and its transpose has one $1$ in each row. Thus, each $b^*_w$
is sent to exactly one element, always the same one for a given equivalence
class.
\qed

For example, the plactic, hypoplactic, and sylvester congruences all satisfy
the statement of Theorem~\ref{th:cong}, so that, we directly conclude that
quotienting, \emph{e.g.}, $\FQSym^*$ by any of these equivalences gives rise
to (already known) Hopf algebras.
Let us first recall some details of these cases before presenting two new
congruences that arose when studying generalizations of permutations.

\subsection{The plactic congruence}

The plactic congruence \cite{LS} is generated by the Knuth relations
\begin{equation}
\left\{ \
\begin{array}{cccl}
aba \, \equiv \, baa & , &
bba \, \equiv \, bab & \qquad \hbox{for} \ a < b  ,  \\
acb \, \equiv \, cab & , &
bca \, \equiv \, bac & \qquad \hbox{for} \ a < b < c \ .
\end{array}
\right.
\end{equation}
The fact that it satisfies the conditions of Theorem \ref{th:cong}
is the starting point of Sch\"utzenberger's proof of the Littlewood-Richardson
rule (see \cite{Loth2} for a simplified presentation), which has been
reinterpreted in \cite{NCSF6} in terms of the algebra $\FSym$ (Free Symmetric
functions), based on standard Young tableaux. 

The plactic monoid remained a mysterious and singular object until its true
nature was revealed with the advent of quantum groups, and in particular of
crystal bases \cite{LLT2,Lit}. 

The number of plactic classes of words of packed evaluation $I$ can be
computed as a sum of Kostka numbers~\cite{Mcd}
\begin{equation}
\sum_{\lambda\vdash n}K_{\lambda\mu}
\end{equation}
if $I$ is a composition of $n$ whose underlying partition is $\mu$.

\subsection{The hypoplactic congruence}

The hypoplactic congruence does for the fundamental basis of quasi-symmetric
functions what the plactic congruence does for Schur functions. It was
originally derived from a degenerate quantum group~\cite{NCSF4}.

It is generated by the plactic relations, and the following
quartic relations (derived from the presentation of the quantum group)
\begin{equation}
\left\{\
\begin{array}{cccl}
baba \, \equiv \, abab & , & baca \, = \, abac & \qquad
\hbox{for} \ a < b < c  ,  \\
cacb \, \equiv \, acbc & , &
cbab \, \equiv \, bacb & \qquad \hbox{for} \ a < b < c  ,  \\
badc \, \equiv \, dbca & , &
acbd \, \equiv \, cdab & \qquad \hbox{for} \ a < b < c < d \ .
\end{array}
\right.
\end{equation}
This may look complicated, but there is a much simpler description: two words
$w$ and $w'$ are congruent if and only if they have the same evaluation, and
their standardizations $\std(w)$ and $\std(w')$ have same recoils~\cite{Nov}.

It is interesting to observe that the hypoplactic congruence is also generated
by the following nonlocal relations
\begin{eqnarray}\label{eq:nlhp}
ac\cdots b &\equiv& ca\cdots b\quad (a\le b<c),\\
b\cdots ac &\equiv& b\cdots ca\quad (a<b\le c),
\end{eqnarray}
where the dots may be any letters~\cite{Nov}.

Hypoplactic classes of words are parametrized by \emph{quasi-ribbon tableaux}.
A \emph{quasi-ribbon tableau} of shape $I$ is a ribbon diagram $r$ of shape
$I$ filled by letters in such a way that each row of $r$ is nondecreasing from
left to right, and each column of $r$ is strictly increasing from \emph{top to
bottom}.
A word is said to be a \emph{quasi-ribbon word} of shape $I$ if it can be
obtained by reading from \emph{bottom to top} and from left to right the
columns of a quasi-ribbon diagram of shape $I$.
Thus, hypoplactic classes of permutations are just recoil classes.
Denoting by $\RC(\sigma)$ the recoil composition of a permutation
$\sigma$ and by $\DC(\sigma)=\RC(\sigma^{-1})$ its descent composition,
we see that the sums
\begin{equation}
R_I = \sum_{\RC(\sigma)=I}\F_\sigma= \sum_{\DC(\sigma)=I}\G_\sigma
\end{equation}
are just the noncommutative ribbon Schur functions~\cite{NCSF1,NCSF6}, and
span a Hopf subalgebra isomorphic to $\Sym$. Dually, the hypoplactic quotient
of $\FQSym$ is isomorphic to $QSym$, the class of a $\F_\sigma$ being
identified with $F_I$, $I=\RC(\sigma)$.

The number of hypoplactic classes of words of packed evaluation $I$
is clearly $2^{\ell(I)-1}$. So hypoplactic classes of packed words are
parametrized by \emph{segmented compostions}, that is, compositions with two
kinds of separators between the parts. There are $3^{n-1}$ segmented
compositions of $n>0$, and the resulting pair of Hopf algebras are the free
tricubical algebra on one generator and its dual~\cite{NTtri}.

Hypoplactic classes of parking functions correspond to \emph{parking
quasi-ribbons}, that is, to quasi-ribbon words which are parking functions.
There is however a simpler parametrization~\cite{NTPark}.

Define a \emph{segmented word} as a finite sequence of non-empty words,
separated by vertical bars, \emph{e.g.}, $232|14|5|746$.

The parking quasi-ribbons can be represented as segmented nondecreasing
parking functions such that the bars can only occur as $\cdots a|b\cdots$ with
$a<b$.

The hypoplactic Hopf algebras obtained from $\PQSym$ have graded dimensions
given by the little Schr\"oder numbers, these are the free triduplicial
algebra on one generator and its dual \cite{NTDup}.

\subsection{The sylvester congruence}

The sylvester congruence \cite{HNT05} is generated by the relations
\begin{equation}
ac\cdots b \equiv ca\cdots b\quad (a\le b<c)
\end{equation}
that is, the first hypoplactic relation in the nonlocal presentation
(\ref{eq:nlhp}).

It is proved in~\cite{HNT05} that the equivalence classes of words under this
congruence are parametrized by the words whose standardized words avoid
the pattern $132$.

The sylvester Hopf algebras obtained from $\FQSym$ are the Loday-Ronco
algebra of planar binary trees, and its (isomorphic) dual.

From $\WQSym$, one obtains the free tridendriform algebra on one generator
(self-dual), whose graded dimensions are given by the little Schr\"oder
numbers (see~\cite{NTPark,NTDup}).

\subsection{Counting sylvester classes}

As for the plactic and hypoplactic congruences, the number of sylvester
classes of words of packed evaluation $I$ can be computed. This relies upon
the noncommutative Lagrange inversion formula \cite{NTLag,NTDup}.

Let $\PBT$ denote the Loday-Ronco algebra of planar binary trees~\cite{LR1},
and $\P_T$ its natural basis (notations are as in~\cite{HNT05}).
Let $\Q_T\in\PBT^*$ be the dual basis of $\P_T$. 
Let $Q_T\in\QSym$ be the image of $\Q_T$ by the adjoint of the Hopf algebra
morphism $\Sym\rightarrow \PBT$ sending $S_n$ to the left comb with $n$
(internal) vertices.
Expanding $Q_T$ in the monomial basis $M_I$ of $\QSym$, we get that
the coefficient of $M_I$ in $Q_T$ is then 1 or 0, according to whether there
is a word of  evaluation $I$ which has sylvester shape $T$ or not. By duality,
the number of terms in the expansion
\begin{equation}
S^I =\sum_{T}\<S^I,\Q_T\> \P_T
\end{equation}
is equal to the number of sylvester classes of evaluation $I$.

Now, 
\begin{equation}
S^I=\sum_{J\le I} R_J,
\end{equation}
and $\<R_J,\Q_T\>$ can be computed from \cite[Prop. 3.2]{NTDup}.
Under the bijection between nondecreasing parking functions and binary trees
described there, each $R_J$ is an interval of the Tamari order, consisting of
all nondecreasing parking functions of packed evaluation $\overline{J}$.
The number of these is the coefficient of $S^{\overline{J}}$ in the Lagrange
series $g$, which is the unique noncommutative symmetric series satisfying the
functional equation
\begin{equation}
g = \sum_{n\ge 0}S_n g^n.
\end{equation}

To summarize, 

\begin{theorem}
The number of sylvester classes of evaluation $I$ is
\begin{equation}
\sum_{J\le\overline{I}}\<M_J,g\>,
\end{equation}
where $\<\cdot ,\cdot \>$ is here the duality bracket between $QSym$ and
$\Sym$.
\end{theorem}
\qed

Defining a basis $E_I$ of $QSym$ by
\begin{equation}
E_I=\sum_{J \le I}M_J
\end{equation}
and denoting by $L_I$ its dual basis in $\Sym$, we can rewrite this number as
\begin{equation}\label{eq:nbslv}
\<E_{\overline{I}},g\>\ =\ \text{coefficient of $L_{\overline{I}}$ in $g$}.
\end{equation}

\begin{example}
{\rm Take $I=(112)$. Then, $\bar I=(211)$, and the term
of degree 4 in $g$ is
\begin{equation}\label{eq:g4}
g_4 = S^{4} + 3S^{31} + 2S^{22} + 3S^{211}
    + S^{13} + 2S^{121} + S^{112} + S^{1111}.
\end{equation}
The compositions $J\le (211)$ are $(211), (31), (22), (4)$. Summing the
coefficients of the corresponding $S^J$, we obtain $3+3+2+1=9$ sylvester
classes of words of evaluation $(112)$. And indeed, representatives of
these classes (avoiding the pattern 132) are
\begin{equation}
1233,\ 2133,\ 2313,\ 2331,\ 3123,\ 3213,\ 3231,\ 3312,\ 3321.
\end{equation}
}
\end{example}

\begin{example}{\rm
The number of sylvester classes of permutations is
\begin{equation}
\<E_{1^n},g\> = \sum_{J\vDash n}\<M_I,g\>
\end{equation}
which is the sum of the coefficients of $g_n$ in the basis
$S^I$, the $n$th Catalan number. Indeed, 
under the specialization $S_n=t^n$ ($A=\{t\}$), the functional equation for $g$
becomes
\begin{equation}
g(t)= \sum_{n\ge 0}t^ng(t)^n = (1-tg(t))^{-1}
\end{equation}
so that $g(t)$ reduces to the generating series of Catalan numbers.
For $n=4$, we can check that $1+3+2+3+1+2+1+1=14$.
}
\end{example}

\begin{example}{\rm
To find the number of sylvester classes of packed words,
we have to compute
\begin{equation}
\left\< \sum_{I\vDash n}E_{\overline{I}},g\right\>=
\left\< \sum_{\lambda\vdash n}2^{n-\ell(\lambda)}m_\lambda,g\right\>
\end{equation}
where the $m_\lambda$ are the monomial symmetric functions.
This can be evaluated by noting that
\begin{equation}
\begin{split}
\sum_{\lambda\vdash n}2^{n-\ell(\lambda)}m_\lambda
&= \frac1{(1-t)^n}
   \sum_{\lambda\vdash n}(1-t)^{\ell(\lambda)}m_\lambda|_{t=1/2} \\
&=\frac1{(1-t)^n} h_n((1-t)X)|_{t=1/2}
\end{split}
\end{equation}
Thus, the number of sylvester classes of packed words
is obtained by putting $x=2$ in the polynomial
\begin{equation}
\begin{split}
N_n(x)
&= \left. \frac1{(1-t)^n}\frac{h_n((1-t)(n+1))}{n+1}\right|_{t=1-1/x} \\
&= \frac1{n+1}
   \sum_{k=0}^n \binom{n+1}{k}\binom{2n-k}{n-k} (1-x)^kx^{n-k}
\end{split}
\end{equation}
a (reciprocal) Narayana polynomial, which gives back the little Schr\"oder
numbers, as expected.
}
\end{example}
\section{New congruences}\label{sec:newcong}

The $m$-generalizations of the classical combinatorial Hopf algebras
introduced in the sequel lead naturally to two new congruences.

\subsection{The hyposylvester congruence}
\label{sec-hypos}

Let $A$ be a totally ordered alphabet.
Let $\sim$ be the equivalence relation on  $A^*$
defined by
\begin{equation}\label{eq:hs1}
w = \cdots ac\cdots  \sim w'=\cdots ca\cdots \ (a<c\in A)
\end{equation}
if a letter $b$ such that $a<b<c$ occurs in $w$, and

\begin{equation}\label{eq:hs2}
w=\cdots ab\cdots a\cdots   \sim w'= \cdots  ba\cdots a\cdots  \ (\ a<b\in A),
\end{equation}
where the dots represent arbitrary letters.
We have
\begin{proposition}
The equivalence relation $\sim$ is a congruence on $A^*$, that is, if
$w\sim w'$, then $uwv\sim uw'v$ for all $u,v\in A^*$. \qed
\end{proposition}

Clearly, $\sim$ is a quotient of the sylvester congruence.
We will call it the \emph{hyposylvester congruence}, and denote it
by $\eqhs$.
Therefore, any hyposylvester class is a disjoint union of sylvester
classes.

The hyposylvester relations are homogeneous, and since they depend only
on \emph{comparisons} between letters, we have:
\begin{lemma}
Two words $w$, $w'$ are in the same hyposylvester class if and only if
\begin{equation}
\ev(w)=\ev(w')
  \quad\text{and}\quad
\pack(w)\equiv_{hs}\pack(w').
\end{equation}
Thus, the number of classes of any evaluation $\alpha$ depends only on the
composition $I=\pack(\alpha)$.
\qed
\end{lemma} 

Let us now give a formula for the number of hyposylvester classes of (packed)
evaluation $I$.

Let us interpret \eqref{eq:hs1} and \eqref{eq:hs2} as rewriting rules,
orienting them from left to right so that $w'>w$ lexicographically.
The number of classes of packed evaluation $I$ is at most the number of words
of evaluation $I$ admitting no such elementary rewriting, since there is at
least one such word in each class. So these maximal words are precisely those
in which each occurence of a letter $j$, except perhaps for the last one, is
followed by smaller letters.
This last occurence is either followed by a smaller letter, or by $j+1$, or by
nothing. Building the word by inserting its letters sequentially in decreasing
order of magnitude, we see that at each step, we have $i_k+1$ choices to place
the letters $k$, since they must necessarily be consecutive, immediately to
the left of an occurence of $k+1$ or at the end of the word.

\begin{example}
{\rm
To find all maximal words of evaluation $I=(222)$, 
we first insert $33$, then a block $22$ in 3 possible ways:
$2233$, $3223$, $3322$, and finally, a block $11$, resulting
into the $3^2$ possibilities
\begin{equation}
112233,\ 211233,\ 223311,\ 311223,\ 321123, \ 322311,\ 331122, \ 332112,\
332211.
\end{equation}
}
\end{example}

Thus, the number of hyposylvester classes of evaluation $I$ is at most
\begin{equation}
\label{eq:numHS}
(i_2+1)(i_3+1)\cdots (i_r+1).
\end{equation}

Let us now prove that there are as many classes as maximal words by providing
a hyposylvester class invariant: the class of any packed word $u$ of
evaluation $I$ is the set of all words $v$ in which the last occurence of each
letter $k$ has exactly as many letters $k+1$ to its left as in $u$. Write
$u\equiv' v$ when this is the case.
Indeed, by a direct check on the elementary rewritings, all words $v$ in the
hyposylvester class of $u$ must have this property, so that $u\eqhs v$ implies
$u\equiv' v$.
But clearly, the number of $\equiv'$-classes of evaluation $I$ is given by
\eqref{eq:numHS}. Thus, $\equiv'$ coincides with $\eqhs$.

Summarizing, we have: 

\begin{theorem}
\label{th:nbhsl}
The number of hyposylvester classes of packed evaluation $(i_1,\ldots,i_r)$
is equal to
\begin{equation}
(i_2+1)(i_2+1)\dots(i_{r}+1).
\end{equation}
The hyposylvester class of a packed word $u$ is the set of all words $v$ in
which the last occurence of each letter $k$ has exactly as many letters $k+1$
to its left as in $u$. Moreover, each class has a unique element that does not
have any elementary rewriting increasing its number of inversions.
\qed
\end{theorem}

\begin{example}{\rm
For permutations, we obtain $2^{n-1}$ classes coinciding with the hypoplactic
classes, since \eqref{eq:hs1}
reduces to the nonlocal hypoplactic relations, 
and \eqref{eq:hs2} is never applied.
}
\end{example}

\begin{example}{\rm
To compute the generating series $c(t)$ of the number $c_n$
of hyposylvester classes of packed words of length $n$,
we can proceed as follows.
Let $\chi_0$ and $\chi_1$ be the characters of $\Sym$ defined by
\begin{equation}\label{eq:defchi}
\chi_0(S_n) = t^n,\quad \chi_1(S_n)=(n+1)t^n\quad (n\ge 0)\,.
\end{equation}
and define a linear map $\alpha$ by
\begin{equation}\label{eq:defalpha}
\alpha(S^{i_1i_2\cdots i_r}) = \chi_0(S_{i_1})\chi_1(S^{i_2\cdots i_r}),
\end{equation}
Let
\begin{equation}
H = \sum_I S^I,\ \text{so that $c(t)=\alpha(H)$}.
\end{equation}
Clearly,
\begin{equation}
H = 1 + \left(\sum_{n\ge 1}S_n\right)H,\ \text{and}\
\alpha(H) = 1 + \chi_0 \left(\sum_{n\ge 1}S_n\right)\chi_1(H)
\end{equation}
so that
\begin{equation}
c(t)= 1+t\cdot\frac1{1-t}\cdot\frac1{1-\sum_{n\ge 1}(n+1)t^n}
=
\frac{1-3t+2t^2}{1-4t+3t^2}.
\end{equation}
The first terms are
\begin{equation}
c(t)=1+t+3\,{t}^{2}+10\,{t}^{3}+34\,{t}^{4}
    +116\,{t}^{5}+396\,{t}^{6}+1352\,{t}^{7}+4616\,{t}^{8}+\dots
\end{equation}
(sequence A007052).
}
\end{example}

\begin{example}{\rm
It will be proved in Section \ref{sec:hypmpark} that
the number of hyposylvester classes of parking functions of length $n$
is equal to the number of ternary trees with $n$ (internal) nodes.
}
\end{example}

As the hyposylvester congruence satisfies the conditions of
Theorem~\ref{th:cong}, there are hyposylvester quotients and subalgebras of
$\FQSym$, $\WQSym$, and $\PQSym$, which will be investigated in the sequel
together with their $m$-analogues. 

\subsection{The metasylvester congruence}

The metasylvester congruence is defined by the nonlocal relations
\begin{eqnarray}
ab\cdots a &\equiv& ba\cdots a\ (a<b),\\
b\cdots ac\cdots b& \equiv& b\cdots ca\cdots b\ (a<b<c)\,.
\end{eqnarray}
These relations are compatible with the concatenation of words,
so that they indeed define a congruence $\eqms$ of the free monoid $A^*$.
The quotient will be called the metasylvester monoid.

This congruence is very close to both the sylvester and hyposylvester
congruences. First, it is also compatible with restriction to intervals
(exchanging $a$ and $c$ requires an extra letter in the interval $[a,c]$ to be
performed), and with any increasing morphism of ordered alphabets (the
relations involve only comparisons between values, not the actual values), so
that it  satisfies the statement of Theorem~\ref{th:cong} as well, and we
can use it to build new Hopf algebras.

We shall use the same technique as for the hyposylvester monoid (see
Section~\ref{sec-hypos}) to provide first a characterization of some canonical
elements of the classes, and then an invariant of the classes, thus proving
all statements of  Theorem~\ref{th:nbms} below in one stroke.

First, consider the packed words $w$ having no elementary rewriting raising
their number of inversions. In other words, consider the words $w$ avoiding
both patterns $ac\dots a$ and $b\dots ac\dots b$.
Then, all the letters $1$ are necessarily consecutive, because of the first
pattern.
Now, if a word $w$ avoids both patterns, its restriction to the letters
greater than $1$ avoids the patterns as well.
From this, we can conclude that the number of metasylvester classes of
packed evaluation $I$ is at most
\begin{equation}
\label{eq-cardms}
(1+i_r) (1+i_{r-1}+i_r) \cdots (1+i_2+\dots+i_r).
\end{equation}

As in the hyposylvester case, we have an invariant of metasylvester classes,
provided by the following algorithm:

\begin{algorithm}
Let $w$ be a word.  The algorithm is defined recursively.

\begin{itemize}
\item Step 1: label the root with the maximal letter $n$ of $w$.
Put $i=n-1$.

\item Step 2: consider the smallest $j>i$ such that the last occurrence of $i$
lies \emph{between} two occurrences of $j$.

\begin{itemize}
\item Step 2.1: if such a $j$ exists, the $(k+1)$th subtree of $j$, where $k$
is the number of occurences of $j$ to the left of the last occurence of $i$
becomes the tree obtained by inserting the restriction of $w$ to $i$ and all
letters that are already in the $(k+1)$th subtree of $j$.
End the algorithm if $i=1$. Otherwise, start again Step 2 with $i=i-1$.

\item Step 2.2: if no such $j$ exists, if this occurrence is to the left of
all maximal letters of $w$, then the leftmost node of the root becomes the
tree obtained by inserting the restriction of $w$ to $i$ and all letters that
are already in the leftmost subtree of the root. The same happens on the
rightmost node if the last occurrence of $i$ is to the right of all
maximal letters of $w$.
End the algorithm if $i=1$. Otherwise, start again Step 2 with $i=i-1$.
\end{itemize}
\end{itemize}
\end{algorithm}

For example, with $w=1413324343$, we obtain Figure \ref{arbredec}.
\begin{figure}[ht]
\begin{equation*}
\vcenter{\xymatrix@C=+6mm@R=+2mm{
*{} & *{} & *{} & {4}\ar@{-}[dlll]\ar@{-}[dl]\ar@{-}[dr]\ar@{-}[drrr] \\
{}  && {1}\ar@{-}[dr]\ar@{-}[dl]\ar@{-}[d] & *{} & {} & *{}
    & {3}\ar@{-}[dll]\ar@{-}[dl]\ar@{-}[d]\ar@{-}[dr]\ar@{-}[drr] \\
*{} & {} & {} & {} & {} & {} & {2}\ar@{-}[dl]\ar@{-}[dr] & {} & {} \\
*{} & *{}& *{} & *{} & *{} & {} & *{} & {} \\
}}
\end{equation*}
\caption{\label{arbredec}The decreasing tree of $1413324343$.}
\end{figure}

Let now be two words $w$ and $w'$ congruent by one elementary
rewriting exchanging a small value $a$ with a greater one. Then all steps
of the algorithm up to $a+1$ are the same, and the
insertion of $a$ follows the same path for both words. Indeed,  $j$ is the
same, or does not exist for both words, and then the induction applies.
So this algorithm provides an invariant, showing that there are at least as
many metasylvester classes of evaluation $I$ as decreasing trees where node
$k$ has $1+i_k$ children, since all such trees can be obtained from at least
one word.

Finally, the number of decreasing trees where node $k$ has $1+i_k$ children
is the same as in Equation~\eqref{eq-cardms} since letter $k-1$ has indeed
$(1+i_k+\dots+i_r)$ places to go in the decreasing tree consisting of nodes
greater than or equal to $k$.

Summarizing, we get:

\begin{theorem}
\label{th:nbms}
The number of metasylvester classes of packed evaluation $(i_1,\ldots,i_r)$
is equal to
\begin{equation}
\label{eq:nbmsl}
(1+i_r) (1+i_{r-1}+i_r) \cdots (1+i_2+\dots+i_r).
\end{equation}
The metasylvester class of a packed word $u$ is the set of all words $v$
giving the same result by the previous algorithm.
Moreover, each class has a unique element admitting no elementary
rewriting increasing its number of inversions.
\end{theorem}

For permutations of size $n$, we get $n!$ and the congruence is trivial.

For packed words of length $n$, the formula yields $\frac12(n+1)!$. Indeed,
one has the induction $P_{n+1}=(n+1)P_n$: write $MS(I)$ for the number of
classes of evaluation $I=(i_1,\dots,i_r)$. Then if $I'=(1+i_1,i_2,\dots,i_r)$
and $I''=(1,I)$ then $MS(I')=MS(I)$ and $MS(I'')=n.MS(I)$. Since all
compositions of size $n+1$ are obtained by one of these two processes, we
conclude to the induction and the formula.

For parking functions of length $n$, the number of classes appears
to be given by Sequence A132624 of ~\cite{Sloane}. It would be interesting to
have a direct proof of this result.

\subsection{Notes on both congruences}

\subsubsection{The canonical elements of the metasylvester congruence}
\label{sec-can-met}
 
Let us first describe two alternative ways to compute the maximal word of a
given metasylvester class.

Consider a decreasing tree where the node labelled $k$ has $1+i_k$ children,
and complete it by adding to each node the missing leaves.  Then read the
nodes following the infix reading: read the leftmost subtree of the root, then
the root, then the next subtree of the root, then the root, and so on, ending
by reading the rightmost subtree of the root.

The word obtained by this algorithm is the same as the following one: starting
from a decreasing $(m+1)$-ary tree, complete it as above, then label each
sector between two edges of a given node by the value of the node and read now
each sector in left-right order.

For example, both readings of the tree represented in Figure~\ref{arbredec}
yield the word $4114433233$. It is indeed the maximal word of its class since
it cannot have any of the patterns $\cdots ac\cdots$ and $b\cdots ac\cdots b$.
We shall call this word the \emph{canonical word} of its class.

\begin{proposition}
The metasylvester canonical words are the words avoiding the pattern $a\cdots
b\cdots a$ with $a<b$ (\emph{i.e.}, a letter cannot be between two occurences
of a smaller one). Such words will be identified with their decreasing trees. 
\end{proposition}

\Proof
Assume first that a word $w$ contains a pattern $a\cdots b\cdots a$.
Choose an occurrence of this pattern such that $a\cdots b$ is of
minimal length $r$.
If $r=2$, this is a pattern $ab\cdots a$, and $w$ cannot be canonical.
If $r>2$, let $x$ be the letter to the left of $b$. Then, $x>a$ or $x<a$.
The first possibility is excluded by the minimality assumption, and
if $x<a$, then $a\cdots xb\cdots a$ is a pattern $b\cdots ac\cdots a$,
so that again, $w$ is not canonical. Conversely, if $w$ does not contain
$a\cdots b\cdots a$, then in particular, it cannot contain the patterns
$ab\cdots a$ or $b\cdots ac \cdots b$.
\qed

\subsubsection{The special case where all $i_k$ are equal}
\label{sec:metaperm}

In the case where all parts $i_k$ of the evaluation are equal to an integer
$m$, which shall be called later the set of $m$-permutations, the resulting
decreasing trees are in fact all decreasing $(m+1)$-ary trees.
In the special cases $m=1$ or $m=2$, one can simplify the output of the
algorithm computing the decreasing tree of a word: one can forget about the
empty subtrees by making the convention that the first subtree is to the left,
the second one to the middle and the third one to the right of their root.
For example, with $m=2$ and $n=9$, one gets the following tree by inserting
$556119367322474898$, whose canonical word is
$559611673322744988$.

\begin{figure}[ht]
\begin{equation*}
\vcenter{\xymatrix@C=+6mm@R=+2mm{
*{}  & *{} & {9}\ar@{-}[dll]\ar@{-}[drr]\ar@{-}[d] \\
{5}  & *{} & {7}\ar@{-}[dr]\ar@{-}[dl]\ar@{-}[d] & *{} & {8} \\
*{}  & {6}\ar@{-}[d] & {3}\ar@{-}[dr] & {4}\\
*{}  & {1} & *{} & {2} \\
}}
\end{equation*}
\caption{\label{arbredect}The decreasing ternary tree of
$556119367322474898$.}
\end{figure}

Note also that in the case of permutations ($m=1$), this algorithm reduces to
the well-known algorithm computing the decreasing (binary) tree of a
permutation.

\subsubsection{Relating the (*)-sylvester congruences}

The hyposylvester classes are unions of sylvester classes, which are
themselves unions of metasylvester classes.
In particular, the canonical elements of the hyposylvester classes are
canonical elements of both the sylvester and metasylvester classes and the
canonical elements of the sylvester classes are canonical elements of the
metasylvester classes.

Hence the hyposylvester classes are in bijection with  a set of decreasing
trees which is  easily  described: these are the decreasing trees such that,
for all $i$ smaller than $n$, node $i$ is a child of $i+1$ but not a rightmost
child, or is the rightmost node of the decreasing tree restricted to nodes
greater than or equal to $i$. This implies clearly that the number of
hyposylvester classes of evaluation $I$ is as in Equation~\eqref{eq:numHS}.

One can do the same with sylvester classes: they are also in bijection with
a set of decreasing trees:  the decreasing trees where, for all  $i\not=n$,
node $i$ cannot be to the left of node $i+1$ and not a child of $i+1$.
One can also understand this result as follows:  there is only one way to
label a given tree so as to obtain a representative of a sylvester class.

\section{Hopf algebras from $m$-Permutations}\label{sec:mperm}

\subsection{Definitions}

We have seen that the Loday-Ronco algebra can be constructed from the
sylvester congruence by summing over sylvester classes of permutations, which
are counted by the Catalan numbers \cite{HNT05}.
Similarly, $(m+1)$-ary trees, which are counted by the $m$-Fuss-Catalan
numbers, are in natural correspondence with sylvester classes of
$m$-permutations.

\begin{definition}
An $m$-permutation of $n$ is a word of evaluation $(m^n)$, that is, a shuffle
of $m$ ordinary permutations of $n$. We denote by $\SG_n^{(m)}$ the set of
$m$-permutations of $n$.
\end{definition}

Their number is therefore $(mn)!/m!^n$.

\begin{lemma}
The map $w\mapsto \std(w)^{-1}$ induces a bijection between $\SG_n^{(m)}$ and
the set of permutations $\tau\in\SG_{mn}$ whose descents are multiples of $m$.
\end{lemma} 

\Proof
Clearly, $\std(w)$ belongs to the shifted shuffle of $n$ factors $12\ldots m$.
\qed

For example, $w=2 2 3 4 1 1 3 4 1 3 4 2$ is a 3-permutation in $\SG_4^{(3)}$.
Its standardization is
$\sigma=4\, 5\, 7\, 10\, 1\, 2\, 8\, 11\, 3\, 9\,  12\, 6$, whose inverse is
\begin{equation}
\sigma^{-1}=
\sruban{5&6&9\\ \ &\ & 1&2&12\\ \ &\ & \ & \ & 3&7&10\\%
        \ & \ &\ & \ & \ & \ &4&8&11\\}
\end{equation}

The restriction of the left weak order to permutations whose descents are
multiples of $m$ induces an order on $m$-permutations which will be called the
(left) $m$-weak order.

\subsection{$m$-Standardization}

Let $A$ be an infinite, totally ordered alphabet. We denote by $\mAn$ the
set of words $w$ of length $mn$ whose number of occurences $|w|_a$ of each
letter is a multiple of $m$.
These words will be called $m$-words.

With an $m$-word $w$, we can associate an $m$-permutation $\std_m(w)$ called
its $m$-standardization. The process is similar to the usual standardization:
scan iteratively the word from left to right, label 1 the first $m$ occurences
of the smallest letter, then 2 the next $m$ occurences of this letter or the
first $m$ occurences of the next one, and so on.
For example, $\std_2(aabacbcbab)=1132535424$.

Define now, for $\alpha\in\SG_n^{(m)}$,
\begin{equation}
\G_\alpha(A)=\sum_{w\in \mAn;\ \std_m(w)=\alpha}w.
\end{equation}
Note that $\G_\alpha$ belongs to $\WQSym$:
\begin{equation}\label{eq:G2M}
\G_\alpha = \sum_{u\in\PW_n^{(m)},\,\std_m(u)=\alpha}\M_u\,.
\end{equation}
where $\PW_n^{(m)}$ denotes the set of packed $m$-words
of length $nm$.

\subsection{The Hopf algebra of Free $m$-quasi-symmetric functions}

Let $\mFQSym_n^*$ be the linear span of the $\G_\alpha$ for
$\alpha\in\SG_n^{(m)}$. 

\begin{proposition}  
The $\G_\alpha$ span a subalgebra of $\Z\<A\>$, and their product 
is given by
\begin{equation}\label{eq:prodG}
\G_\alpha\G_\beta=\sum_{\gamma\in\alpha\star_m\beta}\G_\gamma
\end{equation}
where the $m$-convolution $\alpha \star_m \beta$ is the set of
$m$-permutations $\gamma=uv$ such that $\std_m(u)=\alpha$ and
$\std_m(v)=\beta$.
\end{proposition}

\Proof Clearly, for any word $w$,
\begin{equation}
 \std_m(\pack(w))=\std_m(w),
\end{equation}
so that
\begin{equation}
\G_\alpha\G_\beta
=\sum_{\substack{\std_m(u)=\alpha \\ \std_m(v)=\beta}} \M_u\M_v
=\sum_{\substack{w=uv\\ \std_m(\Park(u))=\alpha\\ \std_m(\Park(v))=\beta}} \M_w
=\sum_{\gamma\in\alpha\star_m\beta}\G_\gamma
\end{equation}
\qed

For example,
\begin{equation}
\G_{2121}\G_{11}= \G_{212133}+\G_{313122}+\G_{323211}.
\end{equation}

\begin{note}
\label{note-prodG}
{\rm
The product $\G_\sigma\G_\tau$ of two $m$-permutations of respective sizes $i$
and $j$ consists in partitioning the set $\{1,\dots,i+j\}$ into two sets $A$
and $B$ of sizes $i$ and $j$ and change $\sigma$ and $\tau$ by increasing
morphisms of alphabets from respectively $\{1,\dots,i\}$ to $A$ and
$\{1,\dots,j\}$ to $B$.
}
\end{note}

Similarly, applying the coproduct of $\WQSym$ yields
\begin{equation}\label{eq:coprodG}
\begin{split}
\Delta\G_\alpha
= \sum_{\std_m(u)=\alpha} \sum_{k=0}^n
       \M_{\pack(u[1..k])}\otimes\M_{\pack(u[k+1..n])}\\
= \sum_{k=0}^n \G_{\std_m(\alpha[1..k])}\otimes\G_{\std_m(\alpha[k+1..n])}
\end{split}
\end{equation}
where for a word $w$, $w[i..j]$ denotes the subword of $w$ whose letters
belong to the interval $[i,j]$.

Summarizing, we have:
\begin{theorem} 
$\mFQSym^*$ is a Hopf subalgebra of $\WQSym$, whose product and coproduct
are given by \eqref{eq:prodG} and \eqref{eq:coprodG}.
\qed
\end{theorem}

\subsection{Commutative image}

The commutative image operation consists in mapping our noncommuting variables
$a_i$ to commuting indeterminates $x_i$.

For $m=1$, the commutative image of $\G_\sigma$ is the quasi-symmetric
function $F_{\RC(\sigma)}$.
There is an analogue of this property for general $m$.
As already observed, for $\alpha\in\SG_n^{(m)}$, the recoils of
$\std(\alpha)$ are all divisible by $m$. Hence, its recoils composition
$I$ is of the form $I=mJ$, where $J$ is a composition of $n$, and we have
\begin{equation}
\G_\alpha(X) = \psi^m(F_J(X))
\end{equation} 
where $\psi^m$ is the algebra morphism $x_i\mapsto x_i^m$ (power-sum plethysm
operator). Indeed, this follows from \eqref{eq:G2M}, since for a packed
$m$-word $u$, $\M_u(X)=\psi^m(M_J)$.

\subsection{Duality}

Let $\F_\alpha$ be the dual basis of $\G_\alpha$  in the dual Hopf algebra
$\mFQSym$.
From \eqref{eq:prodG} and \eqref{eq:coprodG}, we obtain by duality
\begin{equation}
\F_\alpha\F_\beta = \sum_{\gamma\in\alpha\ssh\beta} \F_\gamma,
\end{equation}
and 
\begin{equation}
\label{F-cop-m}
\Delta \F_\alpha
 = \sum_{\substack{\alpha=u\cdot v\\ u,v\in A^{(m)}}} \F_{\std_m(u)} \otimes \F_{\std_m(v)}.
\end{equation}
where $u$ and $v$ run over $m$-word factorizations of $\alpha$.

For example,
\begin{equation}
\begin{split}
\F_{11} \F_{1212} =
\F_{112323} + \F_{121323} + \F_{123123} +
\F_{123213} + \F_{123231}\\
 + \F_{211323} + \F_{213123} + \F_{213213} + \F_{213231} + \F_{231123} \\
+ \F_{231213} + \F_{231231} + \F_{232113} + \F_{232131} + \F_{232311} .
\end{split}
\end{equation}

The calculation of coproducts is easy, since, as $\alpha=u\cdot v$ is an
$m$-permutation, if its prefix $u$ is an $m$-word, then $v$ also is an
$m$-word.
For example, with $m=2$, all $\F_{\alpha}$ with $\alpha$ a permutation of
$1122$ are primitive, except
\begin{equation}
\Delta \F_{1122} = \F_{1122}\otimes1 + \F_{11}\otimes\F_{11}
 + 1\otimes\F_{1122},
\end{equation}
and
\begin{equation}
\Delta \F_{2211} = \F_{2211}\otimes1 + \F_{11}\otimes\F_{11}
 + 1\otimes\F_{2211}.
\end{equation}

\begin{note}
The Hopf algebra is a \slha in the sense of Section~\ref{sub-cong} since its
product is the shifted shuffle and the coproduct amount to an algorithm that
commutes with the action of a group and is compatible with the packing
algorithm.
\end{note}

The product formulas for $\F$ and $\G$ imply that $\mFQSym$ is free and
cofree, so that by~\cite{Foi-free}, it is self-dual.

Note that the number of terms in a product $\F_\alpha\F_\beta$ with
$(\alpha,\beta)$ in $\SG_n^{(m)} \times \SG_p^{(m)}$ is equal to
\begin{equation}
\binom{m(n+p)}{mn}
\end{equation}
whereas the number of terms in a coproduct of an element in $\SG_n^{(m)}$
does not only depend on $m$ and $n$.
On the dual side, for the $\G$ basis, the number of terms in a product or a
coproduct depends only on the sizes of the $m$-permutations.
Thus, an explicit isomorphism between $\mFQSym$ and $\mFQSym^*$ cannot
be expected to be as simple as in the case $m=1$.

\subsection{Noncommutative symmetric functions}

The dual of the Hopf quotient $\psi^m(QSym)$ of $\mFQSym^*$ is the Hopf
subalgebra of $\mFQSym$ freely generated by the
\begin{equation}
S_{n}^{(m)}\mapsto \F_{1^m2^m\ldots n^m}
\end{equation}
whose coproduct is given by
\begin{equation}
\Delta S_{n}^{(m)}=\sum_{i+j=n}S_i^{(m)}\otimes S_j^{(m)}.
\end{equation}
This algebra is clearly isomorphic to $\Sym$, and the ribbon basis is embedded
as
\begin{equation}
R_I^{(m)}=\sum_{C(\std(\alpha)^{-1})=mI}\F_\alpha.
\end{equation}

The adjoint $\phi_m$ of $\psi^m$ is the algebra morphism of $\Sym$ such that
$\phi_m(S_{nm})=S_n$ and $\phi_m(S_k)=0$ if $m$ does not divide $k$. On the
ribbon basis, for a composition $I$ of $mn$,
\begin{equation}
\phi_m(R_I)=\epsilon_m(I)R_K
\end{equation}
where $K$ is the greatest composition of $n$ such that $mK\le I$,
and $\epsilon_m(I)=(-)^{\ell(I)-\ell(K)}$.

For example,
\begin{equation}
\phi_2(R_{23342})=-R_{1321}.
\end{equation}

\subsection{Plactic quotient}

The number of plactic classes of $m$-permutations is given by triangle A188403
of~\cite{Sloane}: Number of $(nm)\times m$ binary arrays with rows in
nonincreasing order, $n$ ones in every column and no more than 2 ones in any
row.
By definition, it is the number of semi-standard tableaux of weight $(m^n)$.
It is also the number of Yamanouchi words of length $mn$ with descents
multiple of $m$, or the number of standard tableaux of size $mn$ with recoils
multiple of $m$.

\subsection{Hypoplactic quotient}

The number of hypoplactic classes of $m$-permutations of $n$ is always
$2^{n-1}$, and the quotients/subalgebras remain isomorphic to $QSym$ and
$\Sym$. 

However, one would have expected to see new algebras in the diagram of
morphisms for $m>1$. These will be provided by the hyposylvester
congruence.

\subsection{Hyposylvester quotient}

Applying Theorem~\ref{th:nbhsl} to evaluation $I=(m^n)$, we obtain:

\begin{corollary}
The number of hyposylvester classes of $m$-permutations of $n$ is 
\begin{equation}
(m+1)^{n-1}.
\end{equation}
\qed
\end{corollary}

For ordinary permutations ($m=1$), we obtain $2^{n-1}$ classes, which in this
case coincide with the hypoplactic classes, since \eqref{eq:hs1} and
\eqref{eq:hs2} imply that two permutations are equivalent iff they have the
same recoils.

\begin{theorem}
The hyposylvester classes of $m$-permutations are intervals of the $m$-weak
order.
\end{theorem}

\Proof
From the above discussion, we know that two $m$-permutations $\alpha$, $\beta$
in the same hyposylvester class have their last occurence of each letter $i$
between the same occurences of $i+1$. The same is true for any element of the
interval $[\alpha,\beta]$.
Moreover, each class has a unique minimal element and a unique maximal
element.
Thus, hyposylvester classes are intervals.
\qed

Since the hyposylvester congruence satisfies the hypotheses of
Theorem~\ref{th:cong}, we obtain Hopf algebras of graded dimensions
$(m+1)^{n-1}$.

The sums of $\F_\alpha$ over hyposylvester classes define
a ribbon-like basis, for which the product is the 
sum of $m+1$ totally associative operations. The resulting algebra
can therefore be identified with the free $As^{(m+1)}$-algebra
on one generator.

Indeed, encode a hyposylvester class as a word $w=w_1\dots w_{n-1}$ on the
alphabet $A_m:=\{0,\dots,m\}$ such that the rightmost letter $i$ has $w_i$
letters $i+1$ to its left.
If an $m$-permutation $\sigma$ satisfies this property, we write $w=\GD(\sigma)$.

Let $w$ be a word on $A_m$ and define
\begin{equation}
\R_{w} := \sum_{\GD(\sigma)=w} \F_\sigma.
\end{equation}
The product $\R_{w} \R_{w'}$ is then a sum of shifted shuffles of $m$-permutations of sizes $p$
and $q$, in which the rightmost letter $p$ will appear with  $k=0,1,\ldots,m$ letters
$p+1$ on its left. Grouping the terms according to the value of $k$ defines
$m+1$ partial products, and we cleary have:

\begin{theorem}
The $\R_w$ generate a Hopf subalgebra $\mNCSF$ of $\mFQSym$.
Their product is given by the formula
\begin{equation}
\R_{w} \R_{w'} := \sum_{k\in A_m} \R_{w\cdot k\cdot w'}.
\end{equation}

Moreover, all operations sending $(\R_w,\R_{w'})$ to
$\R_{w\cdot k\cdot w'}$ are associative and mutually associative, so
that $\mNCSF$ can be identified with the free $As^{(m+1)}$-algebra on one
generator.
\qed
\end{theorem}

Moreover, to mimic the case $m=1$, we shall define $F_w$ as the dual basis of
the $\R_w$ and define the dual of $\mNCSF$ as $\mQSym$. Then consider the set
of all $F_w$ indexed by words on the subalphabet $\{0,m\}$. This set spans
a commutative subalgebra of $\mQSym$: there is a direct isomorphism of Hopf
algebras sending these elements to the usual $F_I$ of $\QSym$. In particular,
this proves that $\mNCSF$ cannot be self-dual since its dual contains a
non-trivial commutative subalgebra. Moreover, a simple computation in the case
$m=2$ also shows that $\mQSym$ is not commutative:
\begin{equation}
F_1 F_\emptyset = F_{10} + F_{02} + F_{21}
\qquad
F_\emptyset F_1 = F_{01} + F_{20} + F_{12}.
\end{equation}

Finally, since the hypoplactic classes are unions of hyposylvester classes,
one obtains an algebra isomorphic to $\QSym$ by quotienting the $F_w$ by the
congruence identifying $0$ with itself and all the other letters with $1$. One
then obtains a binary word encoding the descents of a composition.

\subsection{Sylvester quotient: the Hopf algebra of plane $(m+1)$-ary trees}

\subsubsection{Number of sylvester classes of $m$-permutations}
\label{sub-mndpf}

Applying \eqref{eq:nbslv}, we shall see that the number of sylvester classes
of $m$-permutations in $\SG_n^{(m)}$ is equal to the Fuss-Catalan number
counting the number of $(m+1)$-ary trees with $n$ internal nodes.

Indeed, let us apply this formula to $I=(m^n)$ ($m$-permutations).
Here, $E_{m^n}$ is symmetric, it is
\begin{equation}
E_{m^n}=\sum_{J\le m^n}M_J=\sum_{I\vDash n}M_{mI}=\psi^m(h_n) 
\end{equation}
so that we are reduced to a scalar product in $Sym$:
\begin{equation}
\label{eq:psihg}
\begin{split}
\<\psi^m(h_n), g_{mn}\>
&=\left\< \psi^m(h_n),\frac{1}{mn+1}h_{mn}((mn+1)X)\right\> \\
&=\left\<h_n,\frac{1}{mn+1}h_{n}((mn+1)X)\right\>
\end{split}
\end{equation}
which is
\begin{equation}
\frac{1}{mn+1}h_{n}(mn+1) = \frac{1}{mn+1}\binom{n(m+1)}{n}.
\end{equation}

\subsubsection{Functional equations, parking functions, and $m$-Dyck words}

Denoting by $\phi_m$ the adjoint of $\psi^m$, the number of sylvester
classes of $m$-perm\-utations can be written as well
\begin{equation}
\<\psi^m(h_n), g_{mn}\> = \<h_n,\phi_m(g_{mn})\>.
\end{equation}
Recall that $\psi^m$ is well-defined on $QSym$ ($\psi^m(M_I)=M_{mI}$), and
that $\phi_m$ is defined on $\Sym$, it is the algebra morphism sending
$S_{mn}$ to $S_n$ and the other $S_i$ to $0$.

Let $g^{(m)}=\phi_m(g)\in\Sym$. Then, the coefficient of $S^I$ in $g^{(m)}$ is
the number of nondecreasing parking functions of evaluation $mI$. The
functional equation satisfied by $g^{(m)}$
\begin{equation}
g^{(m)} = \sum_{n\ge 0} S_n \left(g^{(m)}\right)^{mn}
\end{equation}
becomes multiplicity-free if one considers $S_0$ as a letter, and the solution
$X$ of the new functional equation
\begin{equation}
X=S_0 + S_1X^m+S_2 X^{2m}+S_3X^{3m}+\cdots
\end{equation}
is, as in \cite{NTLag},
\begin{equation}\label{eq:gfX}
X = \sum_{\pi\in\NDPF_n^{(m)}}S^{\ev(\pi)\cdot 0^m}
\end{equation}
where $\NDPF_n^{(m)}$ denotes the set of nondecreasing $m$-parking functions
(see Section \ref{sec:mpark} for their definition)
of length $n$, which are in bijection with ordinary nondecreasing parking
functions of length $mn$ whose evaluation is multiple of $m$ (by repeating $m$
times each letter).

For example, with $m=2$
\begin{equation}
X =S^0 + S^{100} + S^{20000}+S^{11000}+S^{10100}+\cdots
\end{equation}

The term of degree $2$ yields the evaluations $(2000)$, $(1100)$,
$(10010)$ corresponding to the $2$-parking functions $11$, $12$, $13$ or to
the ordinary parking functions $1111$ ,$1122$, $1133$.

Introducing the operator 
\begin{equation}
\Omega:\ S^{(i_1,i_2,\ldots,i_r)}\longmapsto S^{(i_1+1,i_2,\ldots,i_r)},
\end{equation}
Equation \eqref{eq:gfX} can be rewritten in the form
\begin{equation}
X = S_0 + (\Omega X)X^m =: S_0 + T_{m+1}(X,X,\ldots,X)
\end{equation}
with an $(m+1)$-linear operator $T_{m+1}$. Thus, the solution by iterated
substitution can be expressed as a sum over $(m+1)$-ary trees, each tree
corresponding to exactly one parking function.

Note that setting $S_n=a^nb$ as in \cite{NTLag}, and $D=Xb$, the functional
equation for $X$ becomes that of the series of $m$-Dyck words
\begin{equation}
Db = b + a(Db)^{m+1}\,.
\end{equation}

\subsubsection{Bijection with $(m+1)$-ary trees}

Now, by definition, sylvester classes of $m$-permutations are parametrized by
the shapes of the decreasing trees of permutations of $mn$ whose descents are
multiple of $m$, whose sum in the basis $\G_\sigma$ is
\begin{equation}
S^{mm\cdots m}=\sum_{I\vDash n}R_{mI}.
\end{equation}
According to \cite{NTDup}, the expansion of a ribbon $R_{mI}$ on the basis
$\P_T$ of $\PBT$ is the sum of the  binary trees corresponding to
nondecreasing parking functions of evaluation $m\bar I$ under the duplicial
bijection of \cite{NTDup}.

Encoding the shape of the binary search tree of an $m$-permutation $\alpha$ by
a nondecreasing parking function and decoding it into an $(m+1)$-ary tree, we
obtain a parametrization of sylvester classes of $m$-permutations by
$(m+1)$-ary trees.
For example, the sylvester class of the 2-permutation $213132$ is encoded by
the binary tree
\setlength\unitlength{2.6mm}
\def\cerp#1#2{\put(#1,#2){\circle*{0.7}}}
\def\cerg#1#2{\put(#1,#2){\circle*{1}}}
\def\filg#1#2#3{\put(#1,#2){\circle*{.7}\Line(#3,#3)}}
\def\film#1#2#3{\put(#1,#2){\circle*{.7}\Line(0,#3)}}
\def\fild#1#2#3{\put(#1,#2){\circle*{.7}\Line(-#3,#3)}}
\def\arbb{{\begin{picture}(1.5,2)
\filg00{1} \fild20{1} \filg{1}{1}{1.5}
\filg{3}01 \fild{4}{1}{1.5}
\cerg{2.5}{2.5}
\end{picture}}}
\def\arbtk{{\begin{picture}(3,4.6)
\film{1.7}0{1.7} \fild{1.7}{1.7}{1.4} \cerg{0.3}{3.1}
\end{picture}}}
\begin{equation}
T = {\arbb}
\end{equation}

Recall from \cite{NTDup} that the duplicial correspondence between binary
trees and nondecreasing parking functions is given by

\begin{equation}
A\backslash B = 
\xymatrix@R=0.1cm@C=2mm{
{\GrTeXBox{A}}\arx1[rd]\\
& {\GrTeXBox{B}}
}
 \rightarrow \P^\alpha\succ\P^\beta,\quad A/B= 
\xymatrix@R=0.1cm@C=2mm{
& {\GrTeXBox{B}}\arx1[ld]\\
{\GrTeXBox{A}} & \\
}
\rightarrow \P^\alpha\prec\P^\beta
\end{equation}
where $\prec$ and $\succ$ are the duplicial operations defined
in~\cite{NTDup}:
\begin{equation}
\P^\alpha\succ \P^\beta=\P^{\alpha\cdot\beta[k]}, \quad
\P^\alpha\prec\P^\beta
=\P^{\alpha\cdot\beta[\max(\alpha)-1]}
\end{equation} 
if $\alpha$ is of length $k$.
For example, $\P^{12}\P^{113}=\P^{12335}$, and
$\P^{12}\prec\P^{113}=\P^{12224}$.

Decomposing $T$ with the over-under operations, we find
\begin{equation}
T = ((\bullet/\bullet)\backslash\bullet)/(\bullet\backslash(\bullet/\bullet))
\end{equation}
so that the corresponding parking function is
\begin{equation}
\tau = ((1\prec 1)\succ 1)\prec (1\succ(1\prec 1)) = 113344.
\end{equation}

The $2$-parking function $134$ has evaluation $(10110)$, and the
corresponding ternary tree is the one evaluating to $S^{10110 00}$ in the
iterative solution of
\eqref{eq:gfX}:
\begin{equation}
\begin{split}
S^{1011000}=T_3(S^0,S^0,S^{1100})
= T_3(S^0,S^0,T_3(S^0,S^{100},S^0))\\
=T_3(S^0,S^0,T_3(S^0,T_3(S^0,S^0,S^0),S^0))
\end{split}
\end{equation}
which is
\begin{equation}
\arbtk
\end{equation}

\subsubsection{$m$-binary trees}

The shapes of the binary trees obtained as the binary search trees of
$m$-permutations are known: these are the left-right flips of the
$m$-binary trees defined in~\cite{these-Pons,CP}.
We shall slightly modify their definition (it amounts to a symmetry through
the vertical axis), so that our
$m$-binary trees are binary trees with the following structure, where all
$T_L$ and $T_R$s are $m$-binary trees.

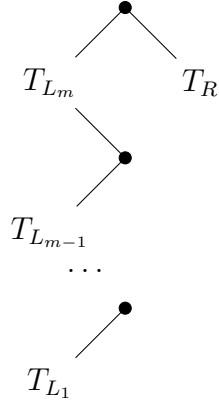
\begin{center}
\tikzstyle{Point} = [fill, radius=0.08]
\begin{tikzpicture}
  \node(x)     at (-1,0){};
  \node(T1)    at (0,-1){$T_R$};
  \node(T2)    at (-2,-1){$T_{L_m}$};
  \node(x1)    at (-1,-2){};
  \node(T3)    at (-2,-3){$T_{L_{m-1}}$};
  \node(inter) at (-1.5,-3.5){$\dots$};
  \node(xm)    at (-1,-4){};
  \node(Tm)    at (-2,-5){$T_{L_{1}}$};

  \draw (x.center) -- (T1);
  \draw (x.center) -- (T2);
  \draw (T2) -- (x1.center);
  \draw (x1.center) -- (T3);
  \draw (xm.center) -- (Tm);
  \draw[Point] (x) circle;
  \draw[Point] (x1) circle;
  \draw[Point] (xm) circle;
\end{tikzpicture}
\end{center}

This definition reveals the recursive structure of these objects, and provides
a direct bijection with $(m+1)$-ary trees : $T_{L_i}$ becomes the $i$-th
subtree of the root and $T_R$ becomes the rightmost (last) subtree of the
root.
If one then restricts the Tamari order to our $m$-binary trees (which form an
upper ideal), one gets the $m$-Tamari lattice defined by Bergeron~\cite{Ber}
up to reversal.

More precisely, denote by $\mR$ the set of permutations whose recoils are
multiple of $m$, so that the number of sylvester classes of $m$-permutations
of size $n$ is equal to the number of sylvester classes of elements of $\mR$
of size $mn$. Note that the elements of $\mR$ form an initial interval of the
weak order, so their sylvester classes also form an initial interval of the
Tamari order. Its greatest element is a very special tree: it is the binary
search tree of the permutation
\begin{equation}
m(n\!-\!1)\!+\!1,\ \dots,\ mn,\
m(n\!-\!2)\!+\!1,\ \dots,\ m(n\!-\!1),\
\dots,\ 1,\ \dots,\ m.
\end{equation}
It is the $(m,n)$ comb-tree, a tree consisting of a root with $n-1$ right
nodes that all have as left subtrees a sequence of $m-1$ left nodes, hence
with $mn$ nodes (see the example for $m=2$ and $n=4$ below).
\begin{figure}[ht]
\begin{equation*}
\vcenter{\xymatrix@C=+7mm@R=+2mm{
*{} &  {}\ar@{-}[dl]\ar@{-}[dr] \\
 {} & *{} &  {}\ar@{-}[dr]\ar@{-}[dl] \\
 *{} &  {} & *{} &  {}\ar@{-}[dl]\ar@{-}[dr] \\
 *{} & *{} &  {} & *{} &  {}\ar@{-}[dl] \\
 *{} & *{} & *{} &  {} & *{} & {} \\
}}
\end{equation*}
\caption{\label{arbremax}The binary search tree of $78563412$.}
\end{figure}

\subsubsection{Hopf algebras and the $m$-Tamari lattice}

Now, since $\mFQSym$ is a \slha and since the sylvester congruence satisfies
the hypothesis of Theorem~\ref{th:cong}, we have

\begin{theorem}
The sums over sylvester classes 
\begin{equation}
\P_T=\sum_{\sylv(\alpha)=T}\F_\alpha
\end{equation}
in $\mFQSym$ span a Hopf algebra $\mPBT$, of graded dimension given by the
Fuss-Catalan sequence. The dual basis $\Q_T$ of $\P_T$ in $\mPBT^*$
is formed of the sylvester classes $\bar\G_\alpha$ for $\sylv(\alpha)=T$.

Moreover, the product of two basis elements $\P_{T'}\P_{T''}$ in $\mPBT$ spans
an interval of the $m$-Tamari order. The bounds of this interval are obtained
through the bijection with $m$-binary trees: order recursively the nodes of
the $(m+1)$-ary trees from left to right in the following way. The first $m$
subtrees in that order, then the root, then the right subtree.
The lower (resp. upper) bound of the product $\P_{T'}\P_{T''}$ consists in
glueing $T$ to the left-most node of $T'$ (resp. $T'$ to the right-most node
of $T$), hence copying the equivalent property in $\PBT$.
\qed
\end{theorem}

\begin{example}{\rm
\setlength\unitlength{1.2mm}
Using the above encoding of sylvester classes by trees, we can now give
a few examples of products and coproduct of ternary trees.
With $m=2$, we have the products
\begin{equation}
\P_{\arbtun} \P_{\arbtun} =
\P_{\arbtdeuxa} + \P_{\arbtdeuxb} + \P_{\arbtdeuxc}.
\end{equation}

\begin{equation}  
\P_{\arbtun} \P_{\arbtdeuxa} =
\P_{\arbta} + \P_{\arbtc} + \P_{\arbtf} + \P_{\arbtb} + \P_{\arbti}.
\end{equation}

\begin{equation}  
\P_{\arbtun} \P_{\arbtdeuxb} =
\P_{\arbte} + \P_{\arbth} + \P_{\arbtj} + \P_{\arbtk}.
\end{equation}

\begin{equation}  
\P_{\arbtun} \P_{\arbtdeuxc} =
\P_{\arbtd} + \P_{\arbtg} + \P_{\arbtl}.
\end{equation}

\begin{equation}  
\P_{\arbtdeuxa}\  \P_{\arbtun} =
\P_{\arbta} + \P_{\arbte} + \P_{\arbtd}.
\end{equation}

\begin{equation}
\P_{\arbtdeuxb} \P_{\arbtun} =
\P_{\arbtc} + \P_{\arbth} + \P_{\arbtg}.
\end{equation}

\begin{equation}
\P_{\arbtdeuxc} \P_{\arbtun} =
  \P_{\arbtf} + \P_{\arbtb}\ + \P_{\arbtj}
+ \P_{\arbti} + \P_{\arbtk} + \P_{\arbtl}.
\end{equation}

\begin{equation}
\P_{\arbtdeuxb} \P_{\arbtdeuxa} =
\P_{\arbqe} + \P_{\arbqd} + \P_{\arbqc} + \P_{\arbqb} + \P_{\arbqa}\,.
\end{equation}

\begin{equation}
\P_{\arbtdeuxa}\ \P_{\arbtdeuxb} =
\P_{\arbqi} + \P_{\arbqh} + \P_{\arbqg} + \P_{\arbqf}.
\end{equation}
and the coproduct
\begin{equation}
\begin{split}
\Delta \P_{\arbqj} &=
\P_{\arbqj}\otimes1 +
(\P_{\arbtdeuxa}+\P_{\arbtdeuxb}+\P_{\arbtdeuxc}) \otimes \P_{\arbtdeuxb} \\
&+ \P_{\arbtun}\otimes (\P_{\arbtb}+\P_{\arbtg}) + 1\otimes\P_{\arbqj}.
\end{split}
\end{equation}

On the $\Q$ basis, we get
\begin{equation}
\Q_{\,\arbtun} \Q_{\,\arbtun} =
\Q_{\,\arbtdeuxa} + \Q_{\,\,\arbtdeuxc}.
\end{equation}

\begin{equation}
\Q_{\,\arbtun} \Q_{\,\arbtdeuxa} =
\Q_{\arbta} + \Q_{\,\arbtf} + \Q_{\,\arbtd}.
\end{equation}

\begin{equation}
\Q_{\,\arbtun} \Q_{\,\,\arbtdeuxb} =
\Q_{\,\arbtb} + \Q_{\,\arbte} + \Q_{\,\,\arbtg}.
\end{equation}

\begin{equation}
\Q_{\,\arbtun} \Q_{\,\,\arbtdeuxc} =
\Q_{\,\arbtd} + \Q_{\,\,\arbti} + \Q_{\,\,\arbtl}.
\end{equation}

\begin{equation}
\Q_{\,\arbtdeuxa} \Q_{\,\arbtun} =
\Q_{\,\arbta} + \Q_{\,\arbtd} + \Q_{\,\,\arbti}.
\end{equation}

\begin{equation}
\Q_{\,\,\arbtdeuxb} \Q_{\,\arbtun} =
\Q_{\,\arbtc} + \Q_{\,\arbtd} + \Q_{\,\arbtk}.
\end{equation}

\begin{equation}
\Q_{\,\,\arbtdeuxc} \Q_{\,\arbtun} =
\Q_{\,\arbtf} + \Q_{\,\arbtd} + \Q_{\,\,\arbtl}.
\end{equation}
}
\end{example}

\begin{note}
\label{mPBT-mults}
{\rm
Finally, as  for $\PBT$, the algebra $\mPBT$ is generated by the
$m$-binary trees whose root has no right child, hence through the bijection,
$\mPBT$ is generated by the $(m+1)$-ary trees with no right (last) subtree.

One can see this directly, since in a product $\P_T \P_{T'}$,
the tree corresponding to the permutation having the maximal number of
inversions is the concatenation of the shifted maximal word of $T'$ and the
maximal word of $T$, hence belongs to the tree where $T'$ is glued to the
rightmost node of $T$ (which is the root when $T$ has no right subtree
attached to it).

The same holds for the same reason for the basis $\Q_T$, except that this time,
the maximal element corresponds to glueing $T$ to the rightmost node of $T'$.
}
\end{note}

\subsubsection{Isomorphism between $\mFQSym$ and a subalgebra of $\FQSym$}

If one standardizes an $m$-permutation on letters $1,\ldots,n$, one obtains a
permutation of size $mn$. The map 
\begin{equation}
\iota_m:\ \F_\alpha\rightarrow \F_{\std(\alpha)} 
\end{equation}
is then an algebra morphism, since the products in the bases $\F$ of $\mFQSym$
and $\FQSym$ are given by the shifted shuffle which is compatible with
standardization.
If one restricts to permutations whose recoil composition has parts
multiple of $m$, we then get an algebra isomorphism.

Note that the coproducts are not the same, so it is only an \emph{algebra}
morphism. It is even worse than that: consider the mapping sending $\F_\sigma$
to itself if its recoils are multiple of $m$ and to $0$ otherwise. Its kernel
is an ideal of $\FQSym$ but not a coideal as one can see by computing the
smallest counter-example
\begin{equation}
\Delta \F_{1423}.
\end{equation}

\subsection{Metasylvester subalgebra and quotient: Hopf algebras of decreasing
$(m+1)$-ary trees}

As observed in Section \ref{sec:metaperm}, metasylvester classes of
$m$-permutations are in bijection with decreasing $(m+1)$-ary trees.
Thus, for $m>1$, new Hopf algebras, based on decreasing $(m+1)$-ary trees can
be constructed. For $m=1$, decreasing binary trees are in bijection with
permutations, and the construction would be only a recoding of $\FQSym$.

\subsubsection{The metasylvester subalgebra}

As already mentioned, we shall identify the canonical metasylvester
$m$-permutations (avoiding the pattern $1\cdots 2\cdots 1$)
with their decreasing trees, which can also be characterized as follows:

\begin{definition}
An $m$-permutation $\tau$ of $n$ is a decreasing $(m+1)$-ary tree if it is
either the empty word, or of the form $u_1nu_2\cdots n u_{m+1}$ where the
packed words of $u_1,\ldots,u_{m+1}$ are decreasing $(m+1)$-ary trees.
We shall denote by $\dT^{(m)}$ the set of these $m$-permutations.
\end{definition}

Up to reversal of the order, this is the same notion as the generalized
Stirling permutations of \cite{JKP}, that correspond to increasing trees.

Regarding $\tau$ as a packed word, one can associate with it a plane tree with
sectors labeled by the letters of $\tau$, as in  \cite{HNT08}. Then, $\tau$
satisfies these conditions iff this tree is a $(m+1)$-ary tree. Labeling the
vertices with the labels of their sectors, one obtains indeed a decreasing
tree.

We can now set, for a decreasing $(m+1)$-ary tree $\tau$
\begin{equation}
\P_\tau = \sum_{\alpha\eqms\tau}\F_\alpha
\end{equation}

For example, with $m=2$, we have
\begin{equation}
\begin{split}
\P_{113322} &= \F_{112332} + \F_{113232} + \F_{113322}. \\
\P_{311322} &= \F_{131232} + \F_{131322} + \F_{311232} + \F_{311322}. \\
\P_{331122} &= \F_{133122} + \F_{131322} + \F_{331122}.
\end{split}
\end{equation}
Note that those three $\P$ elements represent one sylvester class.

As $\eqms$ satisfies the hypotheses of Theorem \ref{th:cong}, we can state:

\begin{theorem}
The $\P_\tau$ span a Hopf subalgebra $\mDT$ of $\mFQSym$.
The dual Hopf algebra $\mDT$ is the metasylvester quotient
${}^m\FQSym^*/\eqms$, based on the equivalence classes
$\Q_\tau =\overline{\G_\alpha}$ ($\alpha\eqms\tau$). 
\end{theorem}

For example, with $m=2$, we have
\begin{equation}
\P_{11} \P_{1122} = \P_{112233} + \P_{211233} + \P_{221133}
                   + \P_{223113} + \P_{223311}.
\end{equation}
\begin{equation}
\P_{11} \P_{2112} = \P_{113223} + \P_{311223} + \P_{321123}
                   + \P_{322113} + \P_{322311}.
\end{equation}
\begin{equation}
\P_{11} \P_{2211} = \P_{113322} + \P_{311322} + \P_{331122}
                   + \P_{332112} + \P_{332211}.
\end{equation}

\begin{equation}
\Delta \P_{311223} = \P_{311223}\otimes 1 + 1\otimes\P_{311223}.
\end{equation}

\begin{equation}
\Delta \P_{322311} = \P_{322311}\otimes 1 + \P_{2112}\otimes \P_{11}
                   + 1\otimes\P_{322311}.
\end{equation}

\begin{equation}
\Delta \P_{113322} = \P_{113322}\otimes 1 + \P_{1122}\otimes\P_{11}
                   + \P_{11}\otimes\P_{2211} + 1\otimes\P_{113322}.
\end{equation}

\begin{equation}
\Q_{11} \Q_{1122} = \Q_{112233} + \Q_{221133} + \Q_{331122}.
\end{equation}
\begin{equation}
\Q_{11} \Q_{2112} = \Q_{113223} + \Q_{223113} + \Q_{332112}.
\end{equation}
\begin{equation}
\Q_{11} \Q_{2211} = \Q_{113322} + \Q_{223311} + \Q_{332211}.
\end{equation}

\begin{equation}
\Delta \Q_{321123} = \Q_{332112}\otimes 1 + \Q_{2112}\otimes\Q_{11}
                   + \Q_{11}\otimes\Q_{2211} + 1\otimes\Q_{321123}.
\end{equation}

As we shall now see, $\mDT$ is also a subalgebra of
${}^m\FQSym^*/\eqms$. This property, which does not hold for the sylvester
congruence, implies immediately the product and coproduct rules for decreasing
trees.

\begin{note}
\label{mDT-mults}
{\rm
As for $\mPBT$ (see Note~\ref{mPBT-mults}), both the algebra
$\mDT$ and its dual $\mDT^*$ are generated by decreasing trees whose roots
have no right (last) subtree. This time, the maximal element of a product
$\P_{\tau} \P_{\tau'}$ corresponds to glueing the decreasing tree $\tau$ to
the right of the shifted decreasing tree $\tau'$.

On the $\Q$ basis, the same property holds and the maximal decreasing tree
corresponds to glueing the decreasing tree $\tau'$ to the right of the shifted
decreasing tree $\tau$.
}
\end{note}

\subsubsection{$\mDT$ as subalgebra of $\mFQSym^*$}

\begin{proposition}
If $\alpha$ and $\beta$ are decreasing $(m+1)$-ary trees, the product
$\G_\alpha\G_\beta$ contains only decreasing $(m+1)$-ary trees. Thus
these trees form a subalgebra of $\mFQSym$.

Moreover, this subalgebra is also a sub-coalgebra.
\end{proposition}

\Proof
The $m$-convolution of two $m$-permutations avoiding the pattern
$1\dots 2\dots 1$ cannot contain this pattern, and the restriction of
such an $m$-permutation to an interval cannot contain it as well.
\qed

For example, we have (with, as always, $m=2$)
\begin{equation}
\G_{11} \G_{2112} = \G_{113223} + \G_{223113} + \G_{332112}.
\end{equation}

\begin{corollary}
(i) The product in $\mDT$ is obtained by selecting the decreasing trees
in the shifted shuffle
\begin{equation} 
\P_{\tau'} \P_{\tau''}
 = \sum_{\tau\in\tau'\ssh\tau'',\ \tau\,\text{tree}}\P_\tau\,.
\end{equation}
(ii)  The coproduct in $\mDT$ is given by 
\begin{equation}
\Delta \P_\tau
 = \sum_{\substack{\tau=u\cdot v\\  \std_m(u), \std_m(v)\in\dT^{(m)}   }}
   \P_{\std_m(u)} \otimes \P_{\std_m(v)}.
\end{equation}
where $u$ and $v$ run over $m$-word factorizations of $\tau$.\\
(iii) The  product in the basis $\Q_\tau$ of $\mDT$ are those
of the basis $\G_\tau$ with $\tau\in\dT^{(m)}$. 
\end{corollary}

\subsubsection{Naked trees}

\begin{proposition}
The shapes of the trees occurring in a product
$\P_{\tau}\P_{\tau'}$
depend only on the shapes of the factors.
The same holds for any coproduct $\Delta \P_\tau$.

Then, the map sending the decreasing trees to their shapes is a Hopf
(quotient) morphism, so that, if one defines $\Q'_T$ as the class of the
decreasing trees $\P_\tau$ of its shape, these $\Q'_T$ span a Hopf quotient
of $\mDT$.
\end{proposition}

\Proof
Let us consider an $(m+1)$-ary tree and consider two decreasing trees of this
shape. Those two trees are obtained from each other by renumbering their
values so that their maximal words are obtained from each other by
multiplication to the left by a permutation. More precisely, all decreasing
posets of the same underlying poset can be obtained by a chain of elementary
transpositions exchanging two consecutive values (see, \emph{e.g.},
\cite{NCSF6} for a refined structure on posets: a $0$-Hecke module
structure).

Thus, we just have to prove that if $\tau$ and $\tau'$ are two
$m$-permutations, both maximal elements of their metasylvester classes, with
the same shape as decreasing trees, and are obtained from one another by an
elementary transposition $\sigma_i$ (acting on their values hence to their
left), then, their left and right products by a given element and their
coproduct when restricted to maximal elements of metasylvester classes happen
to have the same shapes.
Note that it means in particular that the rightmost value $i$ is not between
two values of $i+1$.
Consider the products $\F_\sigma \F_\tau$ and $\F_\sigma \F_{\tau'}$ for
any $\sigma$ maximal in its metasylvester class.
Then, match the elements of both products, according to the positions of the
letters coming from $\sigma$ in the shuffles $\sigma\ssh\tau$ and
$\sigma\ssh\tau'$.
First, either these elements are both maximal elements of metasylvester
classes, or none is. In the first case, they must have the same shape, since
when applying the algorithm computing the decreasing tree of an
$m$-permutation, all steps are identical up to the exchange of $i$ and $i+1$,
because the rightmost $i$ is not between two values of $i+1$.
The same holds for the product to the right by a given $\F_\sigma$.

The discussion of the coproduct is easy.
First, start again with two $m$-permutations, maximal, with the same shape
as decreasing trees, and differing by a single transposition $\sigma_i$.
Then, thanks to Section~\ref{sec-can-met} that gives a way to compute the
maximal element of a class given its decreasing tree, we directly conclude
that those two $m$-permutations can be de-contatenated at the exact same
places.
Moreover, if the last occurrence of $i$ is not between two occurrences of
$i+1$, then it is also the case for any prefix and suffix of these
$m$-permutations, hence proving both coproducts always have same shape.

The Hopf structure is now implied by the previous discussion.
\qed

As an example of the quotient property, we have
\begin{equation}
\P_{223113} \P_{11} =
\P_{22311344} + \P_{22431134} + \P_{42231134} + \P_{42243113}.
\end{equation}

\begin{equation}
\P_{113223} \P_{11} =
\P_{11322344} + \P_{11432234} + \P_{41132234} + \P_{41143223}.
\end{equation}

Note that the quotient basis is defined as $\Q'$ and not as $\P'$. We shall
understand right after the examples why it is so.
Let us now compute a few products of the $\Q'$ with $m=2$.
\setlength\unitlength{1.2mm}
\begin{equation}
\Q'_{\,\arbtun} \Q'_{\arbtdeuxa} =
\Q'_{\arbta} + \Q'_{\arbtc} + \Q'_{\arbtf} + \Q'_{\arbtb} + \Q'_{\arbtd}.
\end{equation}
\begin{equation}
\Q'_{\,\arbtun} \Q'_{\,\arbtdeuxb} =
   \Q'_{\arbtb} + \Q'_{\arbte} + \Q'_{\,\,\arbth}
 + \Q'_{\,\,\arbtj} + \Q'_{\,\,\arbtg}.
\end{equation}
\begin{equation}
\Q'_{\,\arbtun} \Q'_{\,\,\arbtdeuxc} =
   \Q'_{\arbtd} + \Q'_{\,\,\arbtg} + \Q'_{\,\,\arbti}
 + \Q'_{\,\,\arbtk} + \Q'_{\,\,\arbtl}.
\end{equation}

\begin{equation}
\Q'_{\arbtdeuxa} \Q'_{\,\arbtun} =
   \Q'_{\arbta} + \Q'_{\arbtb} + \Q'_{\arbtd}
 + \Q'_{\,\,\arbtg} + \Q'_{\arbte} + \Q'_{\,\,\arbti}.
\end{equation}
\begin{equation}
\Q'_{\,\,\arbtdeuxb} \Q'_{\,\arbtun} =
  \Q'_{\arbtc} + \Q'_{\,\,\arbth} + \Q'_{\,\,\arbtk}.
\end{equation}
\begin{equation}
\Q'_{\,\,\arbtdeuxc} \Q'_{\,\arbtun} =
   \Q'_{\arbtf} + \Q'_{\arbtb} + \Q'_{\arbtd}
 + \Q'_{\,\,\arbtg} + \Q'_{\,\,\arbtj} + \Q'_{\,\,\arbtl}.
\end{equation}

\begin{note}
\label{mPBTp-mults}{\rm
Thanks to Note~\ref{mDT-mults}, it is immediate that the maximal term in a
product $\Q'_T \Q'_{T'}$ is obtained by glueing $T$ as the rightmost child of
the rightmost node of $T'$. So the algebra of the $\Q'$ is free on the trees
whose root do not have a right child.

Note in particular that the maximal term in a product of $\Q'$ is equal to the
maximal term in a product of $\Q$ of $\mPBT$, hence the choice of notations.
Note also that this implies that the algebra of the $\Q'$ is isomorphic
\emph{as an algebra} to $\mPBT^*$.
}
\end{note}

Dually, we have

\begin{theorem}\label{th:DT}
Let $\P'_T$ be the sum of all decreasing trees $\G_\tau$ of shape $T$.
Then, the $\P'_T$ span a Hopf subalgebra of $\mDT^*$ and $\mFQSym^*$.

Moreover, this algebra is isomorphic to $\mPBT$.
\end{theorem}

\Proof
As we have shown before, $\mDT$ is isomorphic to the quotient algebra of
$\mFQSym$ where one sends to zero all the $\F_w$ where $w$ is not the maximal
word of a metasylvester class. So $\mDT$ is isomorphic to the algebra
generated by the $\F_w$ where $w$ are the maximal words of a metasylvester
class, understanding that the products and coproducts are computed by killing
in their expansions the non-maximal words. Moreover, the isomorphism is
direct: send each $\P_\tau$ to $\F_\tau$.
 
Given this, since there is a well-defined quotient on the $\P_\tau$ sending
classes of maximal words to their shape, this means dually that the $\Q'_T$
define a Hopf subalgebra of $\mDT^*$.

Consider now the algebra of the $\P'$, and a product $\P'_T \P'_{T'}$.
Thanks to Note~\ref{mDT-mults}, we see that the maximal term in this product
is obtained by glueing $T'$ as the rightmost child of the rightmost node of
$T'$. So the algebra of the $\P'$ is free on the trees whose root do not have
a right child. Since it is the same for the dual basis $\Q'$, we conclude that
the algebra generated by the $\P'$ is free and cofree of dimension given by
the Fuss-Catalan numbers, and is therefore isomorphic to the algebra $\mPBT$
since there exists at most one free and cofree algebra given its series of
dimensions~\cite{Foi-free}.
\qed

For example, with $m=2$,
\begin{equation}
\P'_{\arbtun} \P'_{\arbtdeuxa} =
\P'_{\arbta} + \P'_{\arbtf} + \P'_{\arbti}.
\end{equation}
\begin{equation}
\P'_{\arbtun} \P'_{\arbtdeuxb} =
\P'_{\arbtb} + \P'_{\arbtk}.
\end{equation}
\begin{equation}
\P'_{\arbtun} \P'_{\arbtdeuxc} =
\P'_{\arbtd} + \P'_{\arbtl}.
\end{equation}

\begin{equation}
\P'_{\arbtdeuxa}\ \P'_{\arbtun} =
\P'_{\arbta} + \P'_{\arbtd}.
\end{equation}
\begin{equation}
\P'_{\arbtdeuxb} \P'_{\arbtun} =
\P'_{\arbtc} + \P'_{\arbtg}.
\end{equation}
\begin{equation}
\P'_{\arbtdeuxc} \P'_{\arbtun} =
\P'_{\arbtf} + \P'_{\arbti} + \P'_{\arbtl}.
\end{equation}

\subsubsection{Functional equations}

In \cite{HNT08}, it has been shown that the basis of the Loday-Ronco
algebra, that is the sums over sylvester classes of permutations, could be
obtained by iterated substitution in a simple functional equation in $\FQSym$:
\begin{equation}\label{eqbinX}
X = 1 + B(X,X)\,,
\end{equation}
where for  $\alpha\in\SG_k$, $\beta\in\SG_l$, and $n=k+l$,
\begin{equation}
B(\G_\alpha,\G_\beta)
 = \sum_{\gf{\gamma=u(n+1)v}{\std(u)=\alpha,\std(v)=\beta}}\G_\gamma\,.
\end{equation}
Clearly,
\begin{equation}
\partial B(\G_\alpha,\G_\beta)=\G_\alpha\G_\beta\,,
\end{equation}
where the derivation $\partial$ is defined by
\begin{equation}
\partial\G_\sigma=\G_{\sigma'},
\end{equation}
$\sigma'$ being the word  obtained by erasing 
the letter $n$ in $\sigma\in\SG_n$. 

The formal sum $X$ of all decreasing $(m+1)$-ary trees in $\mFQSym^*$
satisfies a functional equation generalizing that of \cite{HNT08}
\begin{equation}
X=1+T_{m+1}(X,X,\ldots,X),
\end{equation}
where the $(m+1)$-linear map $T_{m+1}$ is defined, for
$\alpha_1,\ldots,\alpha_{m+1}$ of respective degrees $k_i$, with
$n=k_1+\cdots+k_{m+1}$, by
\begin{equation}
T_{m+1}(\G_{\alpha_1},\G_{\alpha_2},\ldots,\G_{\alpha_{m+1}})
  = \sum_{\pack(u_i)=\alpha_i}\G_{u_1(n+1)u_2(n+1)\cdots (n+1)u_{m+1}}.
\end{equation}
Since $\partial$ is a derivation, the functional
equation lifts the differential equation
\begin{equation}
\frac{dx}{dt} =x^{m+1},\ x(0)=1
\end{equation}
for the exponential generating series
\begin{equation}
x(t)=(1-mt)^{-1/m}
\end{equation}
of the number of decreasing $(m+1)$-ary trees.

The map $\phi: \ \G_\alpha\mapsto t^n/n!$ for $\alpha\in\SG_n^{(m)}$
is a character, and $\phi(X)=x(t)$. Using its $q$-analogue as in \cite{HNT08},
one may easily derive $q$-enumerations of decreasing trees.

The basis $\Q_\tau$ is generated by the $(m+1)$-ary tree solution of the
functional equation.


The sylvester classes of $m$-permutations can also be generated by
a functional equation, as in \cite{HNT08}. Let us illustrate this
for $m=2$ (ternary trees).

Let $X$ and $Y$ be the sums in $\FQSym$ of all $\G_\sigma$ with even descents,
of respectively even and odd lengths.

These series are in $\Sym$:
\begin{eqnarray}
X &=& 1 +S^2+S^{22}+S^{222}\cdots = (1-S_2)^{-1},\\
Y &=& S^1+S^{21}+S^{221}+\cdots = XS_1.
\end{eqnarray}
Note that the number of sylvester classes of evaluation $2^n1$ is
$\frac1{n+1}\binom{3n+1}{n}$ (1, 2, 7 30, 143, $\dots$, which is A006013,
the Hilbert series of the free triplicial algebra on one generator).

With the usual derivation $\partial$ of $\FQSym$ or $\Sym$,
we have the system
\begin{eqnarray}
\partial X &=& YX,\\
\partial Y &=& X + Y^2.
\end{eqnarray}
The fixed point equation is
\begin{eqnarray}
X &=& 1 + B(Y,X),\\
Y &=&  B(Y,Y)+B(X,1).
\end{eqnarray}

Setting $A=t\EE$ in the differential system yields
\begin{eqnarray}
x' &=& yx,\\
y' &=& x + y^2
\end{eqnarray}
with $x(0)=1$ and $y(0)=0$. This implies $(y/x)'=1$ so that 
$y=tx$ and $x=(1-\frac12t^2)^{-1}$ as expected.

Setting $A=t$ in the fixed point system gives back the OGFs for the numbers
of sylvester classes. The system becomes
\begin{eqnarray}
x &=& 1 + tyx,\\
y &=& tx+ty^2 ,
\end{eqnarray}
which implies $x=1+t^2x^3$, the OGF for ternary trees.

In terms of dendriform operations,
\begin{eqnarray}
X &=& 1 + Y\succ \bullet\prec X,\\
Y &=&  Y\succ\bullet\prec Y+X\succ\bullet.
\end{eqnarray}

\section{Hopf algebras from $m$-packed words}\label{sec:mpw}

\subsection{$m$-Packed words}

An $m$-packed word of degree $n$ is  a packed $m$-word, \emph{i.e.},
a packed word of length $mn$ whose evaluation is multiple of $m$.
Their set will be denoted by $\PW^{(m)}_n$.

The noncommutative characteristic of the set of  $m$-packed words
of evaluation $mI$ is $S^{mI}$, so that
\begin{equation}
\ch (\mPW) = \sum_I S^{mI} = (1-S_m-S_{2m}-S_{3m}-\cdots)^{-1},
\end{equation}
and the exponential generating function of $a_n(m)=|\PW^{(m)}_n|$ is
\begin{equation}
\left( 1 - \sum_{n\ge 1}\frac{t^n}{(mn)!}\right)^{-1}
=\sum_{n\ge 0} a_n{(m)}\frac{t^n}{(mn)!}.
\end{equation}

\begin{example}{\rm
For $m=2$, the sequence is \cite[A094088]{Sloane}:
\begin{equation} 	
1, 1, 7, 121, 3907, 202741, 15430207, 1619195761, 224061282907
\end{equation}
whose exponential generating series is
\begin{equation}
\frac{1}{2-\cosh(x)}\quad \text{(even coefficients).}
\end{equation}

We can observe that these numbers are given by evaluating at $x=2$
the reciprocal $2$-Eulerian polynomials, counting $2$-permutations
according to their number of descents (plus one):
\begin{equation}
7 = 5+1\times 2^1 ,\ 121=61+28\times 2^1+2^2, \dots
\end{equation}
Indeed, the natural action of the $0$-Hecke algebra on  $\PW^{(m)}_n$ has
noncommutative characteristic $\sum_{I\vDash n} S^{2I}$, and
\begin{equation}
\sum_I S^{2I}\left(\frac{t}{1-t}\right)^{\ell(I)}
 = \sum_{J}t^{\ell(J)}R_{2J}.
\end{equation} 
}
\end{example}

\subsection{The Hopf algebras $\mWQSym$}

Going back to the definition of $\WQSym$ (see Formulas~\eqref{prodG-wq}
and~\eqref{coprG-wq}), it is immediate that the $\M_u$ for words $u$ whose
evaluation is multiple of $m$ span a Hopf subalgebra of $\WQSym$, which will
be denoted by $\mWQSym$.

Its dual is the quotient of $\WQSym^*$ by the Hopf ideal generated by the
$\NN_u$ such that $u\not\in\PW^{(m)}$.

The sylvester quotient of $\WQSym$ is known to be the free tridendriform
algebra on one generator \cite{NTtri,NT06}.
We shall now investigate the case of $\mWQSym$.
Note that this algebra is a \slha since it is a subalgebra of $\WQSym$ which
is itself a \slha.

\subsection{Sylvester classes}

To find the number of sylvester classes of $m$-packed words, we have to
compute
\begin{equation}
\left\< \sum_{I\vDash n}E_{m\overline{I}},g\right\>=
\left\< \sum_{\lambda\vdash n}2^{n-\ell(\lambda)}m_\lambda,\phi_m(g)\right\>
\end{equation}
One can give directly a $q$-enumeration by remarking that
\begin{equation}
 \sum_{\lambda\vdash n}2^{n-\ell(\lambda)}m_\lambda=
\frac1{(1-t)^n}\sum_{\lambda\vdash n}(1-t)^{\ell(\lambda)}m_\lambda|_{t=1/2}
=\frac1{(1-t)^n}h_n((1-t)X)|_{t=1/2}
\end{equation}
Thus, the number of sylvester classes of $m$-packed words is obtained by
putting $x=2$ in the polynomial
\begin{equation}
\begin{split}
N_n^{(m)}(x)
& = \left. \frac1{(1-t)^n}\frac{h_n((1-t)(mn+1))}{mn+1}\right|_{t=1-1/x}\\
& = \frac1{mn+1}
    \sum_{k=0}^n \binom{mn+1}{k} \binom{(m+1)n-k}{n-k} (1-x)^k x^{n-k}
\end{split}
\end{equation}

For example, with $m=2$, the first polynomials are
\begin{equation}
\begin{split}
& 1,\ \ x+2,\ \ x^2+6 x+5,\ \ x^3+12 x^2+28 x+14, \\
& x^4+20 x^3+90 x^2+120 x+42, \dots
\end{split}
\end{equation}
See Figure~\ref{Nar-2} for a triangular array-like representation of these
polynomials.
Putting $x=2$, we obtain the required numbers: for $m=2$, these are
\begin{equation}
1, 4, 21, 126, 818, 5594, 39693, 289510, \dots
\end{equation}
(sequence A003168, number of ways to dissect a convex $(2n+2)$-gon with
non-crossing diagonals so that no $(2m+1)$-gons ($m>0$) appear).

For $m=3$ we get a new sequence:
\begin{equation}
 1, 5, 34, 267, 2279, 20540, 192350, 1853255
\end{equation}

The coefficients of $N_n^{(m)}$ are the generalized Runyon
numbers~\cite{Cig,Song}.
Their generating function is the functional inverse of
\begin{equation}
f_m(z) = \frac{z}{(1+z)(1+xz)^m}
\end{equation}
Their combinatorial interpretation is known \cite{Cig}: $R_{n,k}^m$ is
$m$-Dyck paths (\emph{i.e.}, paths from $(0,0)$ to $(mn,0)$ by steps $u=(1,1)$
and $d=(1,-m)$) with $k$ peaks.

Note that the reciprocal polynomials $x^nN_n^{(m)}(1/x)$ are the $m$-Narayana
polynomials in the sense of~\cite{Song}.
As we shall see below, their evaluation at $x=2$ yields the number of
hypoplactic classes of $m$-parking functions.

Finally, setting $x=1+q$ instead of just $x=2$ provides the following
$q$-analogues (for $m=2$)
\begin{equation}
\begin{split}
& 1,\ \ q+3,\ \ q^2+8 q+12,\ \ q^3+15 q^2+55 q+55, \\
& q^4+24 q^3+156 q^2+364 q+273, \dots
\end{split}
\end{equation}
See Figure~\ref{Narqb-2} for a triangular array-like representation of these
polynomials.
Their coefficients form the triangle A102537 (number of dissections of a
convex $(2n+2)$-gon by $k-1$ noncrossing diagonals into $(2j+2)$-gons, $1\le
j\le n-1$).

For $m=3$, we obtain new polynomials, encoding a new triangle:
\begin{equation}
\begin{split}
& 1,\ \ q+4,\ \ q^2+11 q+22,\ \ q^3+21 q^2+105 q+140, \\
& q^4+34 q^3+306 q^2+969 q+969, \dots
\end{split}
\end{equation}
See Figure~\ref{Narqb-3} for a triangular array-like representation of these
polynomials.
These are $m$-generalizations of the triangle A033282, obtained for $m=1$
(number of diagonal dissections of a convex $n$-gon into $k+1$ regions). 

\subsection{Sylvester quotient of $\mWQSym$}

Since the algebra $\mWQSym$ is a subalgebra of $\WQSym$, its sylvester
quotient is a subalgebra of the sylvester quotient of $\WQSym$.

Its bases may be labeled by $m$-Schr\"oder paths, which for $m=1$ are in
bijection with reduced plane trees \cite{Song}.

Alternatively, $m$-Schr\"oder paths are in bijection with the natural
generalization of Schr\"oder trees: plane trees in which each internal vertex
has at least $m+1$ children.

\subsection{Hypoplactic quotient of $\mWQSym$}

Since $\ch(\mPW_n)=\sum_{I\vDash n}S^{mI}$, the number of
hypoplactic classes of length $mn$ is
\begin{equation}
\sum_{I\vDash n}2^{\ell(I)-1} = 3^{n-1}.
\end{equation}
The hypoplactic subalgebras of $\mWQSym^*$ and
quotients of $\mWQSym$ are thus always isomorphic to the free tricubical
algebra on one generator and its dual.

\subsection{Hyposylvester quotient of $\mWQSym$}

For $m$-packed words, the generating series of the number of hyposylvester
classes is
\begin{equation}
\begin{split}
c_m(t)
 &= 1 + \sum_{r\ge 1} \sum_{i_1,\ldots,i_r\ge 1}
         (mi_2+1)\cdots(mi_r+1)t^{i_1+\cdots i_r} \\
 &= 1+\frac{t(1-t)}{1-(m+3)t+2t^2}.
\end{split}
\end{equation}
Indeed, with $\alpha$ and $\chi_1$ defined as in \eqref{eq:defalpha}
and \eqref{eq:defchi},
\begin{equation}
c_m(t^m)
 = \alpha(H_m),\ \text{where}\ H_m
 = \sum_I S^{mI}
 = 1+\left(\sum_{n\ge 1}S_{mn}\right) H_m.
\end{equation}
so that
\begin{equation}
c_m(t^m)= 1 + \frac{t^m}{1-t^m}\frac1{1-\sum_{n\ge 1}(mn+1)t^{mn}}\,.
\end{equation}

These series are not in \cite{Sloane}, but their inverses are:
A001835 ($m=2$), A004253 ($m=3$), A001653 ($m=4$),
A049685 ($m=5$), A070997 ($m=6$),
A070998 ($m=7$), A072256 ($m=8$), A078922 ($m=9$),
A077417 ($m=10$), A085260 ($m=11$), A001570 ($m=12$),
A160682 ($m=13$), A157456 ($m=14$), A161595 ($m=15$), \dots

Actually, it can be easily proved that the hyposylvester subalgebras are free,
and these sequences give their numbers of generators by degree (coefficients
of the series $1-c_m(t)^{-1}$).

Thanks to their formula, it is easy to generalize their definition: the
coefficient of $t^n$ in $c_m(t)$ is the number of compositions of $n-1$ with
$m+2$ type of ones and $m$ types of the other values.

\subsection{Metasylvester quotient of $\mWQSym$}

From \eqref{eq:nbmsl}, we get that the number of metasylvester classes
of $m$-packed words in $\PW_n^{(m)}$ is
\begin{equation}
\prod_{k=1}^{n-1}(mk+2).
\end{equation}
The resulting algebras will be investigated in another publication.

\section{$m$-parking functions}

\subsection{Definition}

There are several definitions of generalized parking functions in the
literature.
We shall use the one of Bergeron \cite{Ber}: these are the words over the
positive integers such that the sorted word satisfies $a_i\le k(i-1)+1$.
We shall denote their set by $\PF_n^{(m)}$.
Their number is $(mn+1)^{n-1}$.

\subsection{Frobenius characteristic}

The Frobenius characteristic of the permutation representation of $\SG_n$ on
$\PF_n^{(m)}$ is related to that of ordinary parking functions of length $mn$
whose evaluation is multiple of $m$:
\begin{equation}
\ch(\PF_n^{(m)})
  = \frac1{mn+1}h_n\left((mn+1)X\right)=\phi_m(g_{mn})=\phi_m(g_{mn})\,,
\end{equation}
where $g_{mn}=\ch(PF_{mn})$ (ordinary parking functions) and $\phi_m$ is the
adjoint of the plethysm operator $\psi^m$.
This can be deduced as usual from the noncommutative $0$-Hecke characteristic
$g^{(m)}$, whose functional equation is
\begin{equation}\label{eq:gm}
g^{(m)} = \sum_{n\ge 0}S_n (g^{(m)})^{mn} .
\end{equation}
The solution is given in \cite{NTLag} in the form
\begin{equation}
g^{(m)} = F(r,m)F(r-1,m)^{-1}
\end{equation}
where $r$ is an arbitrary integer, and
\begin{equation}
F(r,m)=\sum_{n\ge 0} S_n((mn+r)A)\,.
\end{equation}
For example,
\begin{equation}
g_3^{(2)} = S_3 + 2S^{21}+4S^{12}+5S^{111}
\end{equation}
the sum of the dimensions of the representations is
\begin{equation}
1+2\times 3 +4\times 3+5\times 6 = 49 = (2\times 3 + 1)^{3-1}
\end{equation}
and the sum of the coefficients is 
\begin{equation}
1+2+4+5 = 12
\end{equation}
the number of nondecreasing $2$-parking functions, see \eqref{eq:2NDPF} below.

As already mentioned, the relation $g^{(m)}=\phi_m(g)$ remains true for the
noncommutative characteristic, where $\phi_m$ is now interpreted as the
adjoint of the Adams operation $\psi^m(M_I)=M_{mI}$ of $QSym$.

\subsection{The Hopf algebra $\mPQSym$}\label{sec:mpark} 

The $m$-parking functions admit Hopf algebra structures. Indeed, there is a
notion of $m$-parkization: to $m$-parkize a word $w$, let $u=w^\uparrow$ be
the sorted word.  If $u_i$ is the first letter violating the constraint
$u_i\le k(i-1)+1$, shift it and all the letter on its right by the minimal
amount $k(i-1)+1-u_i$ so that the constraint gets satisfied, and iterate the
process so as to obtain a nondecreasing $m$-parking function. Finally, put
back each letter at its original place.

Now, let $\a$ be an $m$-parking function and define
\begin{equation}
\G_\a := \sum_{\mPark(w)=\a} w.
\end{equation}
Then, as it is the case for $m=1$ (see Section~\ref{sec-park}), these $\G_\a$
span a subalgebra of $\K\langle A\rangle$.
This algebra is denoted by $\mPQSym^*$.

It is then straightforward that the dual basis $\F_\a$ has its product defined
by shifted shuffle and its coproduct is given by deconcatenation and
$m$-parkization.

In particular, $\mPQSym$ is a \slha on the $\F$ basis.

\subsection{Hyposylvester classes of $m$-parking functions}
\label{sec:hypmpark}

The number of hypo\-syl\-vester classes of $m$-parking functions of size $n$
can be computed as follows.

Let $\chi_0$ and $\chi_1$ be the characters of $\Sym$ defined in
\eqref{eq:defchi}, $\alpha$ as in \eqref{eq:defalpha},
so that the generating series of the number $a_n(m)$ of hyposylvester classes
of $m$-parking functions is
\begin{equation}
X(t) := \sum_{n\ge 0}a_n(m)t^m = \alpha(g_m)\,.
\end{equation}
From \eqref{eq:gm}, we obtain, setting $Y=\chi_1(g_m)$,
\begin{equation}
X(t)=\sum_{n\ge 0}t^n\chi_1(g_m)^{mn} = \frac{1}{1-tY^m},
\end{equation}
since $\chi_1$ is a character. Applying again $\chi_1$ to \eqref{eq:gm},
we obtain a second equation
\begin{equation}
Y = \sum_{n\ge 0}(n+1)t^nY^{mn} = \frac{1}{(1-tY^m)^2}
\end{equation}
and combining both, we have $Y^m = X^{2m}$, so that $X(t)$
satisfies the functional equation
\begin{equation}
X = 1+tX^{2m+1}\quad (X(0)=1).
\end{equation}
Thus, we have proved:

\begin{proposition}
The number of hyposylvester classes of $m$-parking functions is equal to the
number of $(2m+1)$-ary trees with $n$ nodes.
\qed
\end{proposition}

\subsection{Nondecreasing $m$-parking functions}

Ordinary nondecreasing parking functions are counted by Catalan numbers,
and admit simple bijections with noncrossing partitions and binary trees. 
The Hopf algebra~$\CQSym$ obtained by summing parking functions having
the same sorted word can be identified with the free duplicial algebra on one
generator \cite{NTDup}. 

Now, nondecreasing $m$-parking functions of length $n$ are in bijection
with ordinary nondecreasing parking functions of length $nm$ whose evaluation
is multiple of $m$: just repeat each letter $m$ times, \emph{e.g.}, with $m=2$,
$115\rightarrow 111155$.

As already seen in Section \ref{sub-mndpf}, they are enumerated by the
Fuss-Catalan numbers, and are therefore in bijection with $(m+1)$-ary trees or
$m$-Dyck paths.

Arguing as in \cite{NTLag}, we can see that the reciprocal coefficients of
$\phi_m(G)$, where $G$ is the noncommutative $q$-Lagrange series, enumerate 
nondecreasing $m$-parking functions by the shifted sum of their entries.
For example, with $m=2$ and  $n=3$, the enumeration of
nondecreasing $2$-parking functions of size 3
\begin{equation}\label{eq:2NDPF}
111, 112, 113, 114, 115, 122, 123, 124, 125, 133, 134, 135
\end{equation}
by sum of their \emph{entries minus one} is
\begin{equation}
1+q+2q^2+2q^3+3q^4+2q^5+q^6.
\end{equation}
These polynomials do not seem to be a known family of $q$-Fuss-Catalan
numbers.

Nondecreasing $m$-parking functions are also in bijection with $m$-divisible
non-crossing partitions.
Maximal chains in this lattice correspond to $m$-parking functions of length
$n-1$~\cite{Ed1}.

\subsection{Hopf algebras}

The sums $\P^\pi$ of all $m$-parking functions $\F_\a$ with a given reordering
$\pi$ span a Hopf subalgebra $\mCQSym$. Its dual is a commutative algebra
realized as usual (basis $\MM_\pi$). The Adams operation $\psi^m$ map
$\MM_\pi\in\mCQSym^*$ to $\MM_{\pi^m\uparrow}\in\CQSym^*$ (an ordinary
NDPF). The dual map $\phi_m$ is an algebra morphism which maps
$\P^{\pi^m\uparrow}$ to $\P^\pi$ and $\P^\tau$ to 0 if $\tau$ is not
$m$-divisible.

\subsection{Hypoplactic and hyposylvester quotients}

\subsubsection{$m$-Parking functions}

The hypoplactic quotient of $\PQSym$ is a Hopf algebra of graded dimension
given by the little Schr\"oder numbers, not isomorphic to the free
tridendriform algebra on one generator (as it is not self-dual). This algebra,
denoted by $\SQSym$, is the free triduplicial algebra on one generator
\cite{NTDup}.

The number of hypoplactic classes of $m$-parking functions is given by the
Fuss-Schr\"oder numbers obtained as follows.
Let 
\begin{equation}
F(y)=1+\frac{ty^m}{1-ty^m}.
\end{equation}
Then,  
\begin{equation}
P_{n,m}(x) = \frac{1}{mn+1}[t^{mn}]F^{mn+1}
\end{equation}
enumerates nondecreasing $m$-parking functions according to the number of
different letters, and the number of hypoplactic classes is 
\begin{equation}
\frac12 P_{n,m}(2).
\end{equation}
For $m=2$, this is sequence A034015:
\begin{equation}
1,1, 5, 33, 249, 2033, 17485, 156033, 1431281, 13412193, 127840085,
\dots
\end{equation}

The polynomials $P_{n,m}/x$ are
\begin{equation}
\begin{split}
& 1, 2x+1, 5x^2+6x+1, 14x^3+28x^2+12x+1, \\
& 42x^4+120x^3+90x^2+20x+1, \dots
\end{split}
\end{equation}
These are the $m$-Narayana polynomials whose reciprocals $N_n^{(m)}(x)$
evaluated at $x=2$ yield the enumeration of sylvester classes of $m$-packed
words.

They are given by
\begin{equation}
\sum_{i=0}^n \frac1{n+1} \binom{2n+2}{i} \binom{n+1}{i+1} x^i.
\end{equation}
Indeed, nondecreasing $m$-parking functions correspond to $m$-divisible
non-crossing partitions under the standard bijection, and the number of those
is known \cite{Ed1}.

Thus:
\begin{proposition}
The number of hyposylvester classes of $m$-parking functions of length $n$
is
\begin{equation}
\sum_{i=0}^n \frac1{n+1} \binom{2n+2}{i} \binom{n+1}{i+1} 2^i.
\end{equation}
\end{proposition}

\subsection{Narayana triangles}
\subsubsection{Standard $m$-Narayana polynomials}
These triangles show the number of non-decreasing $m$-parking functions with
a given number of different letters (parameter $p$):
Figures~\ref{Nar-1} to~\ref{Nar-3}.

\begin{figure}[ht]
\begin{equation}
\begin{array}{|c|c|c|c|c|c|c|}
\hline
n\backslash p & 1 &  2 &  3 &  4 &  5 & 6 \\
\hline
1   & 1 &&&&&                        \\
\hline
2   & 1 &  1 &&&&                   \\
\hline
3   & 1 &  3 &  1 &&&               \\
\hline
4   & 1 &  6 &  6 &  1 &&          \\
\hline
5   & 1 & 10 & 20 &  6 &  1 &    \\
\hline
6   & 1 & 15 & 50 & 50 & 15 & 1 \\
\hline
\end{array}
\end{equation}
\caption{\label{Nar-1}The usual Narayana triangle.}
\end{figure}

\begin{figure}[ht]
\begin{equation}
\begin{array}{|c|c|c|c|c|c|c|}
\hline
n\backslash p
    & 1 &  2 &  3 &  4 &  5 & 6 \\
\hline
1   & 1     &&&&&                    \\
\hline
2   & 1 &  2     &&&&               \\
\hline
3   & 1 &  6 &  5    &&&           \\
\hline
4   & 1 & 12 & 28 & 14 &&         \\
\hline
5   & 1 & 20 & 90 &120 & 42  &    \\
\hline
6   & 1 & 30 &220 &550 &495 &132\\
\hline
\end{array}
\end{equation}
\caption{\label{Nar-2}The Narayana triangle at $m=2$.}
\end{figure}

\begin{figure}[ht]
\begin{equation}
\begin{array}{|c|c|c|c|c|c|c|}
\hline
n\backslash p
    & 1 &  2 &  3 &   4 &   5 &  6  \\
\hline
1   & 1           &&&&&                  \\ 
\hline
2   & 1 &  3           &&&&             \\
\hline    
3   & 1 &  9 & 12          &&&         \\
\hline
4   & 1 & 18 & 66 &  55       &&      \\
\hline
5   & 1 & 30 &210 & 455 & 273   &    \\
\hline
6   & 1 & 45 &510 &2040 &3060 & 1428\\
\hline
\end{array}
\end{equation}
\caption{\label{Nar-3}The Narayana triangle at $m=3$.}
\end{figure}

\subsubsection{Modified Narayana triangles}

These triangles show the number of hypoplactic classes of
$m$-parking functions according to their number of recoils:
Figures~\ref{Narq-1} to~\ref{Narq-3}.

\begin{figure}[ht]
\begin{equation}
\begin{array}{|c|c|c|c|c|c|c|}
\hline
n\backslash p
    &    1 &    2 &    3 &    4 &    5 &   6 \\
\hline
1   &    1                       &&&&&            \\ 
\hline
2   &    2 &    1                     &&&&       \\
\hline
3   &    5 &    5 &    1                  &&&   \\
\hline
4   &   14 &   21 &    9 &    1              &&\\
\hline
5   &   42 &   84 &   56 &   14 &    1      & \\
\hline
6   &  132 &  330 &  300 &  120 &   20 &   1 \\
\hline
\end{array}
\end{equation}
\caption{\label{Narq-1}The $x=1+q$ Narayana triangle at $m=1$.}
\end{figure}

\begin{figure}[ht]
\begin{equation}
\begin{array}{|c|c|c|c|c|c|c|}
\hline
n\backslash p
    &    1 &    2 &    3 &    4 &    5 &   6 \\
\hline
1   &    1         &&&&&                          \\ 
\hline
2   &    3 &    2       &&&&                     \\
\hline
3   &   12 &   16 &    5    &&&                 \\
\hline
4   &   55 &  110 &   70 &   14  &&              \\
\hline
5   &  273 &  728 &  702 &  288 & 42 &        \\
\hline
6   & 1428 & 4760 & 6160 & 3850 & 1155 & 132 \\
\hline
\end{array}
\end{equation}
\caption{\label{Narq-2}The $x=1+q$ Narayana triangle at $m=2$.}
\end{figure}

\begin{figure}[ht]
\begin{equation}
\begin{array}{|c|c|c|c|c|c|c|}
\hline
n\backslash p
    &    1 &    2 &    3 &    4 &    5 &   6 \\
\hline
1   &    1             &&&&&                      \\ 
\hline
2   &    4 &    3           &&&&                 \\
\hline
3   &   22 &   33 &   12        &&&             \\
\hline
4   &  140 &  315 &  231 &   55    &&          \\
\hline
5   &  969 & 2907 & 3213 & 1547 &  273 &      \\
\hline
6   & 7084 &26565 &39270 &28560 &10200 &1428 \\
\hline
\end{array}
\end{equation}
\caption{\label{Narq-3}The $x=1+q$ Narayana triangle at $m=3$.}
\end{figure}

\subsubsection{Modified Narayana triangles (second version)}

The triangle at $m=1$ is the same as in Figure~\ref{Narq-1} since the Narayana
triangle at $m=1$ is symmetric.
These triangles enumerate the number of sylvester classes of
$m$-packed words by their size minus their number of different letters:
Figures~\ref{Narqb-2} and~\ref{Narqb-3}.

\begin{figure}[ht]
\begin{equation}
\begin{array}{|c|c|c|c|c|c|c|}
\hline
n\backslash p
    &    1 &    2 &    3 &    4 &    5 &   6 \\
\hline
1   &    1      &&&&&                             \\ 
\hline
2   &    3 &    2    &&&&                        \\
\hline
3   &   12 &    8 &    1 &&&                    \\
\hline
4   &   55 &   55 &   15 &    1   &&           \\
\hline
5   &  273 &  364 &  145 &   24 &    1   &    \\
\hline
6   & 1428 & 2380 & 1400 &  350 &   35 &   1 \\
\hline
\end{array}
\end{equation}
\caption{\label{Narqb-2}The reversed $x=1+q$ Narayana triangle at $m=2$.}
\end{figure}

\begin{figure}[ht]
\begin{equation}
\begin{array}{|c|c|c|c|c|c|c|}
\hline
n\backslash p
    &    1 &    2 &    3 &    4 &    5 &   6 \\
\hline
1   &    1                &&&&&                   \\ 
\hline
2   &    4 &    1              &&&&              \\
\hline
3   &   22 &   11 &    1           &&&          \\
\hline
4   &  140 &  105 &   21 &    1       &&       \\
\hline
5   &  969 &  969 &  306 &   34 &    1  &     \\
\hline
6   & 7084 & 8855 & 3850 &  700 &   50 &   1 \\
\hline
\end{array}
\end{equation}
\caption{\label{Narqb-3}The reversed $x=1+q$ Narayana triangle at $m=3$.}
\end{figure}

\section{Conclusions and perspectives}

\subsection{Other algebras}

We have focused our attention on Hopf algebras of trees, and mainly
investigated the sylvester and metasylvester congruences.
A detailed study of the other ones (plactic, hypoplactic, hyposylvester and
stalactic) will be done in another publication. The colored versions, as well
as generalizations of the cambrian congruences need also to be explored.
In the plactic case, one may expect to obtain at some point information on
power-sum plethysms of Schur functions.

\subsection{Bidendriform structures}

All the algebras presented here, as $\mFQSym$, as $\mWQSym$, or as $\mPQSym$
and their sylvester, metasylvester, and hyposylvester subalgebras are
bidendriform bialgebras. This can be proved by means of a generalization of
Theorem~\ref{th:cong} stating that a subalgebra of a bidendriform bialgebra is
itself bidendriform if the monoid satisfies the conditions of
Theorem~\ref{th:cong} and one extra condition: the congruence never
modifies the last letter of a word.
This result and a summary of techniques studying bidendriform structures will
be the subject of a forthcoming paper. 

\subsection{Operads}

Finally, it is highly probable that the new algebras constructed in this paper
are related to some operads, old or new. Their identification would certainly
shed new light on the relationship between operads and combinatorial Hopf
algebras.

\section{Appendix: Tables}


We present here for small $m$ tables of the number of $m$-permutations along
with their numbers of classes for the sylvester, hyposylvester, and
meta-sylvester congruences (Figures~\ref{mPerms} to~\ref{mP-meta}), and the
same data for $m$-packed words (Figures~\ref{mWords} to~\ref{mW-meta}),
$m$-parking functions (Figures~\ref{mPark} to~\ref{mPark-meta}). And, for
comparison, multiple parking functions (shuffles of ordinary parking
functions, Figures~\ref{m-nPark} to~\ref{mnPark-meta}).

When a sequence is in OEIS, its number is written to the right of it (except
for the hyposylvester classes of $m$-permutations enumerated by powers of
integers).


\begin{figure}[ht]
\begin{equation}
\begin{array}{|c||c|c|c|c|c|}
\hline
m\backslash n
    & 1 &  2  &      3 &           4 &              5 \\
\hline
1   & 1 &  2  &      6 &          24 &            120 \\
\hline
2   & 1 &  6  &     90 &        2520 &          113400 \\
\hline
3   & 1 & 20  &   1680 &      369600 &       168168000 \\
\hline
4   & 1 & 70  &  34650 &    63063000 &    305540235000 \\
\hline
5   & 1 & 252 & 756756 & 11732745024 & 623360743125120 \\
\hline
\end{array}
\hskip.5cm
\begin{array}{c}
~ \\
A000142 \\  
A000680 \\  
A014606 \\  
A014608 \\  
A014609 \\  
\end{array}
\end{equation}
\caption{\label{mPerms}$m$-permutations.}
\end{figure}

\pagebreak[4]

\begin{figure}[ht]
\begin{equation}
\begin{array}{|c||c|c|c|c|c|}
\hline
m\backslash n
    & 1 &  2  &   3 &     4 &      5 \\
\hline
1   & 1 &  2  &   5 &    14 &     42 \\
\hline
2   & 1 &  3  &  12 &    55 &    273 \\
\hline
3   & 1 &  4  &  22 &   140 &    969 \\
\hline
4   & 1 &  5  &  35 &   285 &   2530 \\
\hline
5   & 1 &  6  &  51 &   506 &   5481 \\
\hline
\end{array}
\hskip.5cm
\begin{array}{c}
~ \\
A000108 \\  
A001764 \\  
A002293 \\  
A002294 \\
A002295 \\
\end{array}
\end{equation}
\caption{\label{mP-sylv}Sylvester classes of $m$-permutations
                        (Fuss-Catalan numbers).}
\end{figure}

\begin{figure}[ht]
\begin{equation}
\begin{array}{|c||c|c|c|c|c|}
\hline
m\backslash n
    & 1 &  2  &   3 &     4 &      5 \\
\hline
1   & 1 &  2  &   4 &     8 &     16 \\
\hline
2   & 1 &  3  &   9 &    27 &     81 \\
\hline
3   & 1 &  4  &  16 &    64 &    256 \\
\hline
4   & 1 &  5  &  25 &   125 &    625 \\
\hline
5   & 1 &  6  &  36 &   216 &   1296 \\
\hline
\end{array}
\end{equation}
\caption{\label{mP-hypos}Hyposylvester classes of $m$-permutations.}
\end{figure}

\begin{figure}[ht]
\begin{equation}
\begin{array}{|c||c|c|c|c|c|}
\hline
m\backslash n
    & 1 &  2  &   3 &     4 &      5 \\
\hline
1   & 1 &  2  &   6 &    24 &    120 \\
\hline
2   & 1 &  3  &  15 &   105 &    945 \\
\hline
3   & 1 &  4  &  28 &   280 &   3640 \\
\hline
4   & 1 &  5  &  45 &   585 &   9945 \\
\hline
5   & 1 &  6  &  66 &  1056 &  22176 \\
\hline
\end{array}
\hskip.5cm
\begin{array}{c}
~ \\
A000142 \\  
A001147 \\  
A007559 \\  
A007696 \\
A008548 \\
\end{array}
\end{equation}
\caption{\label{mP-meta}Metasylvester classes of $m$-permutations.}
\end{figure}

\pagebreak[3]


\begin{figure}[ht]
\begin{equation}
\begin{array}{|c||c|c|c|c|c|}
\hline
m\backslash n
    & 1 &   2 &      3 &           4 &               5 \\
\hline
1   & 1 &   3 &     13 &          75 &             541 \\
\hline
2   & 1 &   7 &    121 &        3907 &          202741 \\
\hline
3   & 1 &  21 &   1849 &      426405 &       203374081 \\
\hline
4   & 1 &  71 &  35641 &    65782211 &    323213457781 \\
\hline
5   & 1 & 253 & 762763 & 11872636325 & 633287284180541 \\
\hline
\end{array}
\hskip.5cm
\begin{array}{c}
~ \\
A000670 \\  
A094088 \\  %
\\
\\
\\
\end{array}
\end{equation}
\caption{\label{mWords}$m$-packed words.}
\end{figure}

\pagebreak[4]

\begin{figure}[ht]
\begin{equation}
\begin{array}{|c||c|c|c|c|c|}
\hline
m\backslash n
    & 1 &  2  &   3 &     4 &      5 \\
\hline
1   & 1 &  3  &  11 &    45 &    197 \\
\hline
2   & 1 &  4  &  21 &   126 &    818 \\
\hline
3   & 1 &  5  &  34 &   267 &   2279 \\
\hline
4   & 1 &  6  &  50 &   484 &   5105 \\
\hline
5   & 1 &  7  &  69 &   793 &   9946 \\
\hline
\end{array}
\hskip.5cm
\begin{array}{c}
~ \\
A001003 \\  
A003168 \\  
\\
\\
\\
\end{array}
\end{equation}
\caption{\label{mW-sylv}Sylvester classes of $m$-packed words.}
\end{figure}

\begin{figure}[ht]
\begin{equation}
\begin{array}{|c||c|c|c|c|c|}
\hline
m\backslash n
    & 1 &  2  &   3 &     4 &      5 \\
\hline
1   & 1 &  3  &  10 &    34 &    116 \\
\hline
2   & 1 &  4  &  18 &    82 &    374 \\
\hline
3   & 1 &  5  &  28 &   158 &    892 \\
\hline
4   & 1 &  6  &  40 &   268 &   1796 \\
\hline
5   & 1 &  7  &  54 &   418 &   3236 \\
\hline
\end{array}
\hskip.5cm
\begin{array}{c}
~ \\
A007052 \\  
A052913 \\  
\\ 
\\
\\
\end{array}
\end{equation}
\caption{\label{mW-hypo}Hyposylvester classes of $m$-packed words}
\end{figure}

\begin{figure}[ht]
\begin{equation}
\begin{array}{|c||c|c|c|c|c|}
\hline
m\backslash n
    & 1 &  2  &   3 &     4 &      5 \\
\hline
1   & 1 &  3  &  12 &    60 &    360 \\
\hline
2   & 1 &  4  &  24 &   192 &   1920 \\
\hline
3   & 1 &  5  &  40 &   440 &   6160 \\
\hline
4   & 1 &  6  &  60 &   840 &  15120 \\
\hline
5   & 1 &  7  &  84 &  1428 &  31416 \\
\hline
\end{array}
\hskip.5cm
\begin{array}{c}
~ \\
A001710 \\  
A002866 \\  
A034000 \\  
A000407 \\  
A034323 \\  
\end{array}
\end{equation}
\caption{\label{mW-meta}Metasylvester classes of $m$-packed words}
\end{figure}



\begin{figure}[ht]
\begin{equation}
\begin{array}{|c||c|c|c|c|c|}
\hline
m\backslash n
    & 1 &  2  &   3 &     4 &      5 \\
\hline
1   & 1 &  3  &  16 &   125 &   1296 \\
\hline
2   & 1 &  5  &  49 &   729 &  14641 \\
\hline
3   & 1 &  7  & 100 &  2197 &  65536 \\
\hline
4   & 1 &  9  & 169 &  4913 & 194481 \\
\hline
5   & 1 & 11  & 256 &  9261 & 456976 \\
\hline
\end{array}
\hskip.5cm
\begin{array}{c}
~ \\
A000272 \\  %
A052750 \\  
A052752 \\  %
A052774 \\  %
A052782 \\  %
\end{array}
\end{equation}
\caption{\label{mPark}$m$-parking functions}
\end{figure}

\begin{figure}[ht]
\begin{equation}
\begin{array}{|c||c|c|c|c|c|}
\hline
m\backslash n
    & 1 &  2  &   3 &     4 &      5 \\
\hline
1   & 1 &  3  &  13 &    69 &    417 \\
\hline
2   & 1 &  5  &  40 &   407 &   4797 \\
\hline
3   & 1 &  7  &  82 &  1239 &  21810 \\
\hline
4   & 1 &  9  & 139 &  2789 &  65375 \\
\hline
5   & 1 & 11  & 211 &  5281 & 154661 \\
\hline
\end{array}
\end{equation}
\caption{\label{mPark-sylv}Sylvester classes of $m$-parking functions
                           (none in~\cite{Sloane}).}
\end{figure}

\begin{figure}[ht]
\begin{equation}
\begin{array}{|c||c|c|c|c|c|}
\hline
m\backslash n
    & 1 &  2  &   3 &     4 &      5 \\
\hline
1   & 1 &  3  &  12 &    55 &    273 \\
\hline
2   & 1 &  5  &  35 &   285 &   2530 \\
\hline
3   & 1 &  7  &  70 &   819 &  10472 \\
\hline
4   & 1 &  9  & 117 &  1785 &  29799 \\
\hline
5   & 1 & 11  & 176 &  3311 &  68211 \\
\hline
\end{array}
\hskip.5cm
\begin{array}{c}
~ \\
A001764 \\  
A002294 \\  
A002296 \\
A059967 \\
A230388 \\
\end{array}
\end{equation}
\caption{\label{mPark-hypo}Hyposylvester classes of $m$-parking functions.}
\end{figure}

\begin{figure}[ht]
\begin{equation}
\begin{array}{|c||c|c|c|c|c|}
\hline
m\backslash n
    & 1 &  2  &   3 &     4 &      5 \\
\hline
1   & 1 &  3  &  11 &    45 &    197 \\
\hline
2   & 1 &  5  &  33 &   249 &   2033 \\
\hline
3   & 1 &  7  &  67 &   741 &   8909 \\
\hline
4   & 1 &  9  & 113 &  1649 &  26225 \\
\hline
5   & 1 & 11  & 171 &  3101 &  61381 \\
\hline
\end{array}
\hskip.5cm
\begin{array}{c}
~ \\
A001003 \\  
A034015 \\  
\\
\\
\\
\end{array}
\end{equation}
\caption{\label{mPark-hypop}Hypoplactic classes of $m$-parking functions.}
\end{figure}

\pagebreak[4]

\begin{figure}[ht]
\begin{equation}
\begin{array}{|c||c|c|c|c|c|}
\hline
m\backslash n
    & 1 &  2  &   3 &     4 &      5 \\
\hline
1   & 1 &  3  &  14 &    87 &    669 \\
\hline
2   & 1 &  5  &  45 &   585 &   9944 \\
\hline
3   & 1 &  7  &  94 &  1879 &  50006 \\
\hline
4   & 1 &  9  & 161 &  4353 & 158035 \\
\hline
5   & 1 & 11  & 246 &  8391 & 386211 \\
\hline
\end{array}
\hskip.5cm
\begin{array}{c}
~ \\
A132624 \\  
\\  %
\\
\\
\\
\end{array}
\end{equation}
\caption{\label{mPark-meta}Metasylvester classes of $m$-parking functions (the
first row in conjecturally A132624).}
\end{figure}


\pagebreak[4]


\pagebreak[4]

\begin{figure}[ht]
\begin{equation}
\begin{array}{|c||c|c|c|c|c|}
\hline
m\backslash n
    & 1 &   2 &      3 &           4 &               5 \\
\hline
1   & 1 &   3 &     16 &         125 &            1296 \\
\hline
2   & 1 &   7 &    136 &        5293 &          347776 \\
\hline
3   & 1 &  21 &   1933 &      483209 &       257484501 \\
\hline
4   & 1 &  71 &  36136 &    68501421 &    349901224576 \\
\hline
5   & 1 & 253 & 765766 & 12012527625 & 648203695298171 \\
\hline
\end{array}
\hskip.5cm
\begin{array}{c}
~ \\
A000272 \\  
\\  %
\\
\\
\\
\end{array}
\end{equation}
\caption{\label{m-nPark}Multiparking functions (shuffles of $m$ ordinary PFs).}
\end{figure}

\begin{figure}[ht]
\begin{equation}
\begin{array}{|c||c|c|c|c|c|}
\hline
m\backslash n
    & 1 &  2  &   3 &     4 &      5 \\
\hline
1   & 1 &  3  &  13 &    69 &    417 \\
\hline
2   & 1 &  4  &  24 &   179 &   1532 \\
\hline
3   & 1 &  5  &  38 &   360 &   3919 \\
\hline
4   & 1 &  6  &  55 &   628 &   8235 \\
\hline
5   & 1 &  7  &  75 &   999 &  15262 \\
\hline
\end{array}
\end{equation}
\caption{\label{mnPark-sylv}Sylvester classes of multiparking functions
(none in~\cite{Sloane}, not even the first one since $n=6$ gives 2759).}
\end{figure}

\begin{figure}[ht]
\begin{equation}
\begin{array}{|c||c|c|c|c|c|}
\hline
m\backslash n
    & 1 &  2  &   3 &     4 &      5 \\
\hline
1   & 1 &  3  &  12 &    75 &    273 \\
\hline
2   & 1 &  4  &  21 &   126 &    818 \\
\hline
3   & 1 &  5  &  32 &   233 &   1833 \\
\hline
4   & 1 &  6  &  45 &   382 &   3498 \\
\hline
5   & 1 &  7  &  60 &   579 &   6017 \\
\hline
\end{array}
\hskip.5cm
\begin{array}{c}
~ \\
A001764 \\  
A003168 \\  
\\
\\
\\
\end{array}
\end{equation}
\caption{\label{mnPark-hypo}
         Hyposylvester classes of multiparking functions.}
\end{figure}

\begin{figure}[ht]
\begin{equation}
\begin{array}{|c||c|c|c|c|c|}
\hline
m\backslash n
    & 1 &  2  &   3 &     4 &      5 \\
\hline
1   & 1 &  3  &  14 &    87 &    669 \\
\hline
2   & 1 &  4  &  27 &   254 &   3048 \\
\hline
3   & 1 &  5  &  44 &   551 &   8919 \\
\hline
4   & 1 &  6  &  65 &  1014 &  20598 \\
\hline
5   & 1 &  7  &  90 &  1679 &  40977 \\
\hline
\end{array}
\hskip.5cm
\begin{array}{c}
~ \\
A132624 \\  
\\  %
\\
\\

\end{array}
\end{equation}
\caption{\label{mnPark-meta}
         Metasylvester classes of multiparking functions.}
\end{figure}


\end{document}